\documentclass[a4paper,8pt]{article}
\usepackage{graphicx}
\usepackage{amsfonts}

\usepackage{caption}
\usepackage{subcaption}
\usepackage{amsmath}
\usepackage{amssymb}
\usepackage{physics}
\usepackage[utf8]{inputenc}
\usepackage[english]{babel}
 \usepackage{textcomp}
 \usepackage{multirow}
 \usepackage{mathrsfs,amsmath}
\usepackage[square,numbers]{natbib}
\newcommand{\dif}{\mathop{}\!\mathrm{d}}
\newcommand{\supp}{\mathop{}\!\mathrm{supp}}
\newcommand{\eup}{\mathop{}\!\mathrm{e}}

\title{Bibliography management: \texttt{natbib} package}
\author{Share\LaTeX}
\date { }

       \author{Yifan Du}
       \title{Turbulence Generation from a stochastic wavelet model}
       
\begin{document}
     
       \maketitle
        \begin{abstract}
      This research presents a new turbulence generation method based on stochastic wavelets and tests its various properties in both homogeneous and inhomogeneous turbulence. Turbulence field can be generated with less basis compared to previous synthetic Fourier methods. Adaptive generation of inhomogeneous turbulence is achieved by scale reduction algorithm and lead to smaller computation cost. The generated turbulence shows good agreement with input data and theoretical results.
      \end{abstract}
       \section{Introduction}
       Turbulence generation has been an important research topic in fluid mechanics area for decades. Before wide application of Direct Numerical Simulation(DNS) and Large Eddy Simulation(LES) in fundamental turbulence research, explicit synthesis scheme was constructed to study turbulence related phenomenon. \citet{kraichnan1970diffusion} proposed a divergence-free synthesis method using random wavelets and applied it on diffusion and particle dispersion in Homogeneous Isotropic Turbulence(HIT). Subsequent improvement and modification soon turned it into a viable tool for both theoretical research and acoustic related computation(\citet{juves1999stochastic}, \citet{fung1992kinematic}). With the increase of computational power, DNS and LES became realistic and successfully enabled people to acquire knowledge from more complex flow phenomenon(\citet{rogallo1984numerical}, \citet{moin1998direct}). Meanwhile, turbulence generation techniques was gradually modified to generate high-fidelity inlet boundary conditions for DNS and LES, as well as interface generation in Detached Eddy Simulation (\citet{spalart2009detached}), rather than directly used for turbulent flow prediction. Synthetic Random Fourier Method(SRFM) has become one of the most important methods for inflow turbulence generation.   
       However, despite its wide application in turbulence related research, SRFM suffers several drawbacks. Fourier basis is naturally not proper for representation of inhomogeneous turbulence because of its global properties. Although \citet{le1997direct} proposed a transform to map isotropic homogeneous turbulence generated from SRFM to general turbulence field with any given Reynolds stress, the transformed turbulence field is no longer incompressible(\citet{wu2017inflow}). Also, SRFM uses global Fourier basis for all wavenumbers, which lead to very large computational cost. 
       
In the meantime, wavelet noise tool was gradually developed for fluid simulation in Computer Graphics area and movie industry. \citet{perlin1985image} constructed the widely used turbulence function that could describe band-limited noise. The form of Perlin Noise was very close to a series of random wavelets. \citet{cook2005wavelet} rigorously constructed band-limited wavelet noise upsampling and downsampling proccedure. \citet{bridson2007curl} used Perlin Noise as vector potential to generate incompressible flow field. \citet{kim2008wavelet} constructed high-resolution incompressible flow field based on the wavelet noise of \citet{cook2005wavelet}. As a new mathematical tool developed in 1980s, wavelet exhibits many delicate and fine properties. Similar to Fourier series, wavelet series can be complete basis for $\mathbb{L}^2$ function space. Unlike Fourier transform which can only extract global frequency(wavenumber) information from functions, wavelet transform can contain both frequency(wavenumber) and location information, which is an appropriate tool for inhomogeneous turbulence(\citet{farge1992wavelet}). 

In this research, a turbulence generation method based on stochastic wavelets is developed and tested in both homogeneous and inhomogeneous cases. Turbulence generation method is proposed in Sec.\ref{Methodology}, including generation system, boundary condition treatment, input spectrum and scale reduction algorithm. Numerical results are presented in Sec.\ref{Numerical Results} for both isotropic homogeneous turbulence and fully developed channel flow to validate its accuracy. An comprehensive discussion and analysis of the proposed method is made in Sec.\ref{Conclusion and Discussion}. An proof of Reynolds Stress preservation of this method is presented in Sec.\ref{Appendix}.        
       
       \section{Methodology}\label{Methodology}
       \subsection{Turbulence Generation Process}

       \subsubsection{Generation of static turbulence field}
        The model presented in this work allows to generate inhomogeneous stochastic turbulent velocity field based on RANS data. Many previous research(\citet{shur2014synthetic}, \citet{lund1998generation}) suggested that the fluctuation velocity field could be expressed in the following form:
      \begin{equation}\label{eq:1}
       \textrm{\textit{\textbf{u}}}= \textrm{\textit{\textbf{A}}}  \cdot(\nabla \times \textrm{\textit{\textbf{M}}})
      \end{equation}
      Where $\textrm{\textit{\textbf{u}}}=(u, v, w)$ is the fluctuation velocity. $\textrm{\textit{\textbf{u}}}$ is constructed in the manner of equation Eq.\ref{eq:1} so that corresponding second order momentum $\langle \textrm{\textit{\textbf{u}}} \textrm{\textit{\textbf{u}}}\rangle$ is equal to Reynolds stress tensor.  $\textrm{\textit{\textbf{A}}}(x)$ is Cholesky decomposition of Reynolds stress tensor $\textrm{\textit{\textbf{R}}}(x)$(\citet{JARRIN2009435}, \citet{shur2014synthetic}) :
       \begin{equation}\label{eq:2}
       \textrm{\textit{\textbf{R}}} = \textrm{\textit{\textbf{A}}}^{T}\textrm{\textit{\textbf{A}}}
       \end{equation}
       $\textrm{\textit{\textbf{A}}}^{T}$ denotes the transpose of $\textrm{\textit{\textbf{A}}}$. $ \textrm{\textit{\textbf{M}}}={(M_x, M_y, M_z)}$ is vector potential field.  $\nabla \times \textrm{\textit{\textbf{M}}}$ term is constructed to be divergence free as suggested in \citet{shur2014synthetic}. In \citet{kim2008wavelet}, vector potential field was constructed using wavelet noise in \citet{cook2005wavelet}. Similarly $\textrm{\textit{\textbf{M}}}$ is decomposed into a sum of wavelet modes:
       \begin{equation}\label{eq:3}
       \textrm{\textit{\textbf{M}}}(\textrm{\textit{\textbf{x}}})=
       \sum_{\vert \textrm{\textit{\textbf{k}}} \vert \in K}
       \sum^{N_i}_{\textrm{\textit{\textbf{x}}}_p}
       q_{\textrm{\textit{\textbf{x}}}_p,\textrm{\textit{\textbf{k}}}}\textrm{\textit{\textbf{O}}}_{\textrm{\textit{\textbf{x}}}_p,\textrm{\textit{\textbf{k}}}}(\mathbf{\omega}_{\textrm{\textit{\textbf{x}}}_p,\textrm{\textit{\textbf{k}}}} \Psi_{\textrm{\textit{\textbf{k}}}}(\textrm{\textit{\textbf{x}}}-\textrm{\textit{\textbf{x}}}_p))
       \end{equation}
       $\textrm{\textit{\textbf{x}}}_p=(x_p,y_p,z_p)$ is the position of wavelet basis in physical space. $\textrm{\textit{\textbf{k}}} = (k_x, k_y, k_z)$ is wavenumber of wavelet basis. Definition of such wavenumber is according to \citet{perrier1995wavelet}. $K = \{ l_1, l_2, ..., l_M \}$ is a series of magnitude of wavenumber vector. For each $l_i$, $\textrm{\textit{\textbf{k}}}$ is randomly chosen on a sphere of radius $l_i$ in spectral space. Such a construction coincides with the nature of energy spectrum function, i.e. $E(l)$ is turbulent kinetic energy distributed on a sphere of radius $l$ in spectral space.  $\textrm{\textit{\textbf{x}}}_p$ is the location of $N_i$ randomly distributed wavelets. The number $N_i$ for each wavenumber $l_i$ is determined using the following expression deducted from uniformly distributed wavelets:  
       \begin{equation}\label{eq:4}
       N_i = [\vert \Omega \vert \big (\frac{l_i}{k_0} \big )^3]
       \end{equation}
Where $\Omega$ represents the flow domain. $\mathbf{\omega}_{\textrm{\textit{\textbf{x}}}_p,b}=(\omega^1_{\textrm{\textit{\textbf{x}}}_p,b},\omega^2_{\textrm{\textit{\textbf{x}}}_p,b},\omega^3_{\textrm{\textit{\textbf{x}}}_p,b})$ is a random vector series assumed to be normally distributed with the following statistics:
      \begin{equation}\label{eq:5}
      \langle \omega^i_{\textrm{\textit{\textbf{x}}}_p,\textrm{\textit{\textbf{k}}}}\rangle = 0
      \end{equation}
      \begin{equation}\label{eq:6}
      \langle \omega^i_{\textrm{\textit{\textbf{x}}}_p,\textrm{\textit{\textbf{k}}}}\omega^j_{\mathbf{x}_p,\textrm{\textit{\textbf{k}}}}\rangle = \delta_{ij}
      \end{equation}
       $\Psi(\textrm{\textit{\textbf{x}}})$ is 3D wavelet basis function constructed in the following tensor-product way:
       \begin{equation}\label{eq:7}
       \Psi_{\textrm{\textit{\textbf{k}}}}(\textrm{\textit{\textbf{x}}}-\textrm{\textit{\textbf{x}}}_p) = \psi (\frac{k_x}{k_0}(x-x_p))\psi (\frac{k_y}{k_0}(y-y_p))\psi (\frac{k_z}{k_0}(z-z_p))
       \end{equation}
       From the analytical result in \citet{deriaz2005towards}, $\Psi_{\textrm{\textit{\textbf{k}}}}(\textrm{\textit{\textbf{x}}}-\textrm{\textit{\textbf{x}}}_p)$ in the above form might not be a complete basis of $\mathbb{L}^2(\mathbb{R}^3)$. However, the basis function above is chosen because of its localization in both physical and spectral space, which offers an appropriate tool for description and synthesis of turbulence  (\citet{farge1992wavelet}). Wavelet function of enough high order cancellation is chosen as 1-D wavelet function $\psi(\cdot)$ as stated in \citet{farge1992wavelet} and \citet{perrier1995wavelet}. $k_0$ is the Fourier wavenumber of wavelet function where its Fourier spectrum reaches peak, i.e.:
       \begin{equation} \label{eq:8}
       \mathscr{F}\{\psi(x)\}(k_0) = \max_{k \in \mathbb{R}}{\mathscr{F}\{\psi(x)\}(k)}
       \end{equation}
       Wavelet function $\psi(x)$ is localized in both physical and spectral space, which represents a local structure of turbulent field with certain band width of wavenumber. $k_0$ defined in Eq.\ref{eq:8} characterizes the most energetic wavenumber of such structure. 
       $\textrm{\textit{\textbf{O}}}$ in Eq.\ref{eq:3} is a random rotation matrix in 3-dimensional space, i.e. $\textrm{\textit{\textbf{O}}} \in SO(3)$.  An efficient way of generating uniformly distributed random rotation is from \citet{stuelpnagel1964parametrization}, where random rotation matrix is generated from random quaternion: 
       $$
       \textrm{\textit{\textbf{O}}}
       =\left( \begin{array}{ccc}
1-2c_u^2-2d_u^2 & 2b_uc_u-2a_ud_u& 2b_ud_u+2a_uc_u \\
2b_uc_u+2a_ud_u & 1-2b_u^2-2d_u^2 & 2c_ud_u-2a_ub_u \\
2b_ud_u-2a_uc_u & 2c_ud_u+2a_ub_u & 1-2b_u^2-2c_u^2
\end{array} \right)
       $$
       where $a_u$, $b_u$, $c_u$, $d_u$ are components of a unit quaterion:
       $$
       a_u = \frac{a}{\vert \textrm{\textit{\textbf{q}}} \vert}, b_u = \frac{b}{\vert \textrm{\textit{\textbf{q}}} \vert},  c_u = \frac{c}{\vert \textrm{\textit{\textbf{q}}} \vert}, d_u = \frac{d}{\vert \textrm{\textit{\textbf{q}}} \vert}
       $$
       where $\textrm{\textit{\textbf{q}}} = a + b\textrm{\textit{\textbf{i}}} + c\textrm{\textit{\textbf{j}}} + d\textrm{\textit{\textbf{k}}}$ is a random quaterion with $a,b,c,d \sim N(0,1)$. 
        
        $q_{\textrm{\textit{\textbf{x}}}_p,\textrm{\textit{\textbf{k}}}}$ is a series of normalized weights to maintain local spectrum property of RANS data:
       \begin{equation}\label{eq:9}
       q_{\textrm{\textit{\textbf{x}}}_p,\textrm{\textit{\textbf{k}}}} = \sqrt{\frac{E(l)\Delta l}{2k_t c_{l}}}
      \end{equation}
       In which $E(l)$ represents local energy spectrum, which is a known function used as input. Various of different spectrum can be used as input spectrum to characterize multiscale feature of turbulence field. $\Delta l$ is difference between two neighboring wavenumber magnitudes in $K$ in Eq.\ref{eq:3}. $k_t$ is the turbulent kinetic energy with the following relation:
       \begin{equation}\label{eq:10}
       k_t = \int^{\infty}_{0}E(l) \dif l 
       \end{equation}
       $c_l$ is a coefficient determined using Monte Carlo method:
       \begin{equation}\label{eq:11}
       c_l = \frac{N_i}{\vert \Omega \vert} \langle \int_{S}\frac{\partial}{\partial x}\Psi_{\textrm{\textit{\textbf{k}}}}\dif V \rangle = \frac{N_i}{\vert \Omega \vert} \langle \int_{S}\frac{\partial}{\partial y}\Psi_{\textrm{\textit{\textbf{k}}}}\dif V \rangle = \frac{N_i}{\vert \Omega \vert} \langle \int_{S}\frac{\partial}{\partial z}\Psi_{\textrm{\textit{\textbf{k}}}}\dif V \rangle
      \end{equation}
       where the ensemble average$\langle \cdot \rangle$ is performed on random variable $\textrm{\textit{\textbf{k}}}$. Such ensemble average diverges without restriction on $\textrm{\textit{\textbf{k}}}$. In real flow problems, flow domain with finite size can only contain wavelet modes with finite large support, which prevents any component of $\textrm{\textit{\textbf{k}}}$ approaching zero and removes singularity in Eq.\ref{eq:11}. $S = \supp(\Psi_{\textrm{\textit{\textbf{k}}}})$ is the support set of tensor-product basis function. It should be noted that depending on the choice of wavelet basis $\psi(\cdot)$, $\Psi_{\textrm{\textit{\textbf{k}}}}$ may not be compactly supported, in which case S is the effective support set of function $\Psi_{\textrm{\textit{\textbf{k}}}}$ in numerical sense, which could be defined as follow:
       \begin{equation}\label{eq:12}
       S = \{ (\textrm{\textit{\textbf{x}}},t) \vert  |\Psi_{\textrm{\textit{\textbf{k}}}}(\textrm{\textit{\textbf{x}}},t)|< \delta \}
       \end{equation}
       where $\delta$ is a small positive number. 
       
       Construction of such system including Eq.\ref{eq:1}, \ref{eq:2}, \ref{eq:3}, \ref{eq:4}, \ref{eq:7}, \ref{eq:9} coincides with the multi-scale and inhomogeneous nature of turbulence, which needs further explanation. It has been stated clearly in multiple literatures that wavelets are good tools for performing energy decomposition to find possible atoms in physical-spectral space(\citet{farge1992wavelet}, \citet{farge2001coherent}). In Eq.\ref{eq:3} velocity vector potential 
$\textrm{\textit{\textbf{M}}}$ is decomposed into a series of wavelet basis $\Psi_{\textrm{\textit{\textbf{k}}}}(\textrm{\textit{\textbf{x}}}-\textrm{\textit{\textbf{x}}}_p)$ with its own characteristic wavenumber $\textrm{\textit{\textbf{k}}}$ and position $\textrm{\textit{\textbf{x}}}_p$. Each wavelet mode represents a vortex structure localized both around position $\textrm{\textit{\textbf{x}}}_p$ in physical space and wavenumber $\textrm{\textit{\textbf{k}}}$ in spectral space. Each $\Psi_{\textrm{\textit{\textbf{k}}}}(\textrm{\textit{\textbf{x}}}-\textrm{\textit{\textbf{x}}}_p)$ is equipped with a random rotation matrix $\textrm{\textit{\textbf{O}}}_{\textrm{\textit{\textbf{x}}}_p,\textrm{\textit{\textbf{k}}}}$ to make sure Eq.\ref{eq:3} is invariant under rotation. Summation on position $\textrm{\textit{\textbf{x}}}_p$ indicates a layer of vortices with same magnitude of characteristic wavenumber at different position. The whole fluctuation velocity field is a superposition of layers of local structures with different magnitudes of wavenumbers. The preservation of Reynolds Stress tensor in the construction of Eq.\ref{eq:1}, \ref{eq:2}, \ref{eq:3} is in Appendix

Such construction of fluctuation velocity field resembles the Synthetic Eddy Methods(SEM) used to generate inlet flow conditions of LES in previous researchs(\citet{JARRIN2009435}, \citet{poletto2013new}, etc.). However, unlike wavelet functions used in this research, the spectral and physical space properties of functions used to generate structures of different scales  in SEM remain unknown. Also, numbers of modes with different scales used in eddy synthesis often need to be determined by experience and tests. Eq.\ref{eq:4} gives a quantitative representation of number of wavelet modes in this system, which comes from the density of wavelets used to completely cover each scale in wavelet theory(\citet{hernandez1996first}, \citet{deriaz2005towards}).  Intuitively speaking, Eq.\ref{eq:4} indicates that turbulence field with larger size need more wavelet modes to cover. Also, there are much more small scale structure with higher wavenumber than large scales structures with lower wavenumber. 

Many previous research have constructed multi-scale system of turbulence using Fourier basis(\citet{fung1992kinematic}, \citet{juves1999stochastic}), and successfully simulated isotropic homogeneous turbulence from it. However, in general cases of anisotropic inhomogeneous turbulence, such construction encounters serious problems. Fourier basis is global in the whole turbulence field, thus the construction procedure in \citet{fung1992kinematic} is not applicable for generation of turbulence with anisotropy and inhomogeneity. \citet{billson2003jet}, \citet{shur2014synthetic} modified Fourier based reconstruction system so that it could be used to generate inhomogeneous turbulence by multiplying each Fourier basis with a weighing function which quantifies the distribution of turbulent kinetic energy on each scale locally. However, such modification changes the properties of Fourier basis and causes aliasing between different wavenumbers. In this research, such problem is solved by using random wavelets which are local in both spectral space and physical space rather than random Fourier basis.

           \subsubsection{Generation of dynamic turbulence field}\label{Generation of dynamic turbulence field}
     For high Reynolds number turbulence, structures of large and small scales might behave differently according to their own kinematic and dynamic properties, thus need to be dealt with differently(\citet{pope2001turbulent}, \citet{fung1992kinematic}, \citet{lafitte2014turbulence}). Large scale contains most of turbulent kinetic energy of the whole turbulence field, while small scale include inertial subrange and viscous subrange. The separation of large and small scales can be achieved by introducing a cutoff wavenumber $k_c$. However, previous research on turbulence generation indicates that such $k_c$ might be difficult to determine a priori. Here $k_c$ is determined with the following relation:
\begin{equation}\label{eq:13}
k_c = 2k_e
\end{equation}
where $k_e$ is the wavenumber where maximum energy spectrum occurs. The magnitudes of wavenumber in $K$ in Eq.\ref{eq:3} can be separated into large scales and small scales. For large scales, time advance is achieved by  advection induced by mean velocity field $\textrm{\textit{\textbf{U}}}$ computed from RANS model and a random Gaussian advective velocity $W_{\textrm{\textit{\textbf{k}}}}$(\citet{fung1992kinematic}, \citet{poletto2013new}):
\begin{equation}\label{eq:14}
\textrm{\textit{\textbf{x}}}_{p}(t) = \textrm{\textit{\textbf{x}}}_{0} + \int_0^{t'}  \textrm{\textit{\textbf{U}}}(t')+ \textrm{\textit{\textbf{W}}}_{\textrm{\textit{\textbf{k}}}}(t') \dif t'
\end{equation}
where $\textrm{\textit{\textbf{x}}}_0$ is initial position of a large scale mode. It is reasonable to assume that $\textrm{\textit{\textbf{W}}}_{\textrm{\textit{\textbf{k}}}}$ has zero mean and variance equal to that of velocity field of mode ${\textrm{\textit{\textbf{k}}}}$:
\begin{equation}\label{eq:15}
\langle\textrm{\textit{\textbf{W}}}_{i,\textrm{\textit{\textbf{k}}}}^2 \rangle = \frac{E(l)\Delta l}{2k_t}\langle u_i u_i \rangle
\end{equation}
where $l=\vert \textrm{\textit{\textbf{k}}} \vert$, $i=1,2,3$ represents three space coordinates and does not imply summation on repeated index. 

For small scales structures, \citet{fung1992kinematic} and \citet{lafitte2014turbulence} suggest that small scales vortices are advected by large scale structures and mean velocity field. A similar formula as Eq.\ref{eq:15} can be written as follow:
\begin{equation}\label{eq:16}
\textrm{\textit{\textbf{x}}}_{p}(t) = \textrm{\textit{\textbf{x}}}_{0} + \int_0^{t'}  \textrm{\textit{\textbf{U}}}(t')+ \textrm{\textit{\textbf{u}}}_{l}(t') \dif t'
\end{equation}
where $\textrm{\textit{\textbf{u}}}_{l}$ means large scale velocity field. For small scale basis, Eq.\ref{eq:1}, \ref{eq:2}, \ref{eq:3} still hold. However, to predict right Lagrangian decorrelation process, basis function $\Psi_{\textrm{\textit{\textbf{x}}}_p,\textrm{\textit{\textbf{k}}}}$ is modified to the following shape:

 \begin{equation}\label{eq:17}
       \Psi_{\textrm{\textit{\textbf{k}}}}(\textrm{\textit{\textbf{x}}}-\textrm{\textit{\textbf{x}}}_p,t) = \psi (\frac{k_x}{k_0}(x-x_p))\psi (\frac{k_y}{k_0}(y-y_p))\psi (\frac{k_z}{k_0}(z-z_p)) \psi(\frac{\omega_{\textrm{\textit{\textbf{k}}}}}{k_0} t - \phi_p)
       \end{equation}
       where $\phi_p$ is a random phase, randomly distributed in support set of function $\psi(\frac{\omega_{\textrm{\textit{\textbf{k}}}}}{k_0} t)$ with uniform distribution. $\omega_{\textrm{\textit{\textbf{k}}}}$ is a time frequency related to wavenumber $\textrm{\textit{\textbf{k}}}$. Generally structure with small wavenumber varies in a slow frequency and vice versa. In inertial subrange, it is assumed that energy at each wavenumber $\textrm{\textit{\textbf{k}}}$ is spread over a range of frequency around a characteristic frequency related to characteristic wavenumber. Thus here for each characteristic wavenumber $\textrm{\textit{\textbf{k}}}$ a random frequency $\omega_\textrm{\textit{\textbf{k}}}$ is generated with the following distribution:
       \begin{equation}\label{eq:18}
       f_{\omega} = \frac{1}{\sqrt{2\pi}\sigma_{\omega}(\textrm{\textit{\textbf{k}}})}\eup^{\frac{-(\omega-\varpi(\textrm{\textit{\textbf{k}}}))^2}{2\sigma^2_{\omega}(\textrm{\textit{\textbf{k}}})}}
       \end{equation}
where $\sigma_{\omega}(\textrm{\textit{\textbf{k}}})=\varpi(\textrm{\textit{\textbf{k}}})=\epsilon^{\frac{1}{3}}\vert \textrm{\textit{\textbf{k}}} \vert^{\frac{2}{3}}$ are variance and  mean of random frequency $\omega_\textrm{\textit{\textbf{k}}}$(\citet{leslie1975developments})

\subsubsection{Incompressible condition modification}
The system Eq.\ref{eq:1}, \ref{eq:2}, \ref{eq:3}, \ref{eq:9} does not satisfy incompressible condition as the result of transformation Eq.\ref{eq:1}. Previous Synthetic Random Fourier Methods(SRFM) all suffer this problem. Some modified SRFM could unify divergence-free constraint, anisotropy and inhomogeneity, but they often involve some special parameter, which largely undermines their generality(\citet{wu2017inflow}, \citet{smirnov2001random}, \citet{huang2010general}, \citet{castro2013time}, \citet{yu2014fully}). A slight modification of formulas allows to generate inhomogeneous turbulence field which satisfies divergence-free constraint. The modified formulas are as follow:

\begin{equation}\label{eq:19}
\textrm{\textit{\textbf{u}}}= \nabla \times \textrm{\textit{\textbf{M}}}
\end{equation}
\begin{equation}\label{eq:20}
       \textrm{\textit{\textbf{M}}}(\textrm{\textit{\textbf{x}}})=
       \sum_{\vert \textrm{\textit{\textbf{k}}} \vert \in K}
       \sum^{N_i}_{\textrm{\textit{\textbf{x}}}_p}
       q_{\textrm{\textit{\textbf{x}}}_p,\textrm{\textit{\textbf{k}}}}\textrm{\textit{\textbf{O}}}_{\textrm{\textit{\textbf{x}}}_p,\textrm{\textit{\textbf{k}}}}(\mathbf{\omega}_{\textrm{\textit{\textbf{x}}}_p,\textrm{\textit{\textbf{k}}}} \Psi_{\textrm{\textit{\textbf{k}}}}(\textrm{\textit{\textbf{x}}}-\textrm{\textit{\textbf{x}}}_p))
       \end{equation}
\begin{equation}\label{eq:21}
       q_{\textrm{\textit{\textbf{x}}}_p,\textrm{\textit{\textbf{k}}}} = \sqrt{\frac{E(l)\Delta l}{2 c_{l}}}
      \end{equation}
The system Eq.\ref{eq:19}, \ref{eq:20}, \ref{eq:21} is incompressible. Also, it considers the inhomogeneity of turbulence field. However, this construction does not preserve local Reynolds stress, only preserves turbulent kinetic energy distribution.

\subsection{Boundary Conditions}
\begin{figure}[!h]\label{fig:1}
\centering
\begin{subfigure}{.5\textwidth}
\centering
\includegraphics[width=1.3\linewidth]{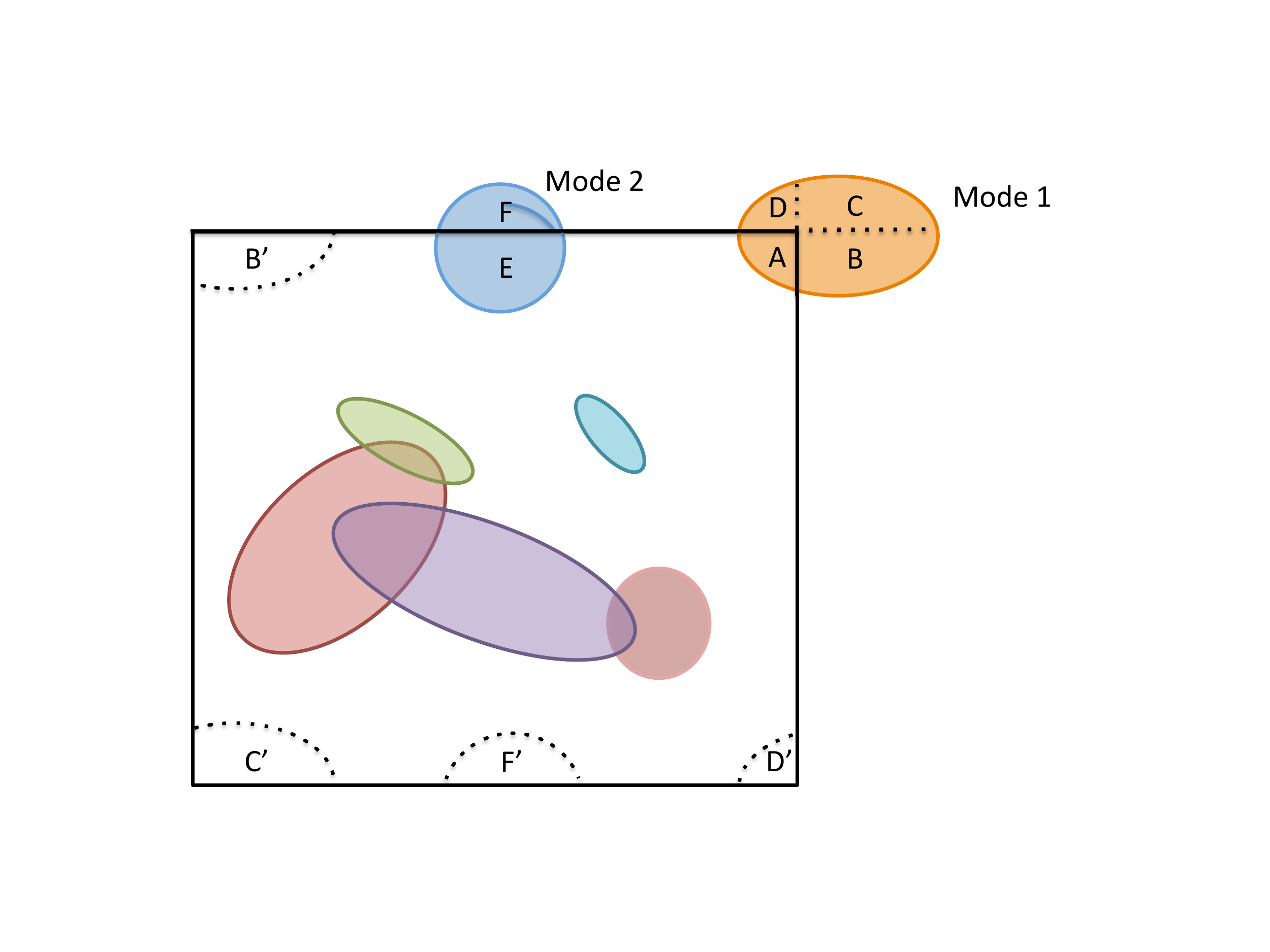}
\caption{Periodic Boundary}
\label{fig:1sub1}
\end{subfigure}%
\begin{subfigure}{.5\textwidth}
\centering
\includegraphics[width=1.3\linewidth]{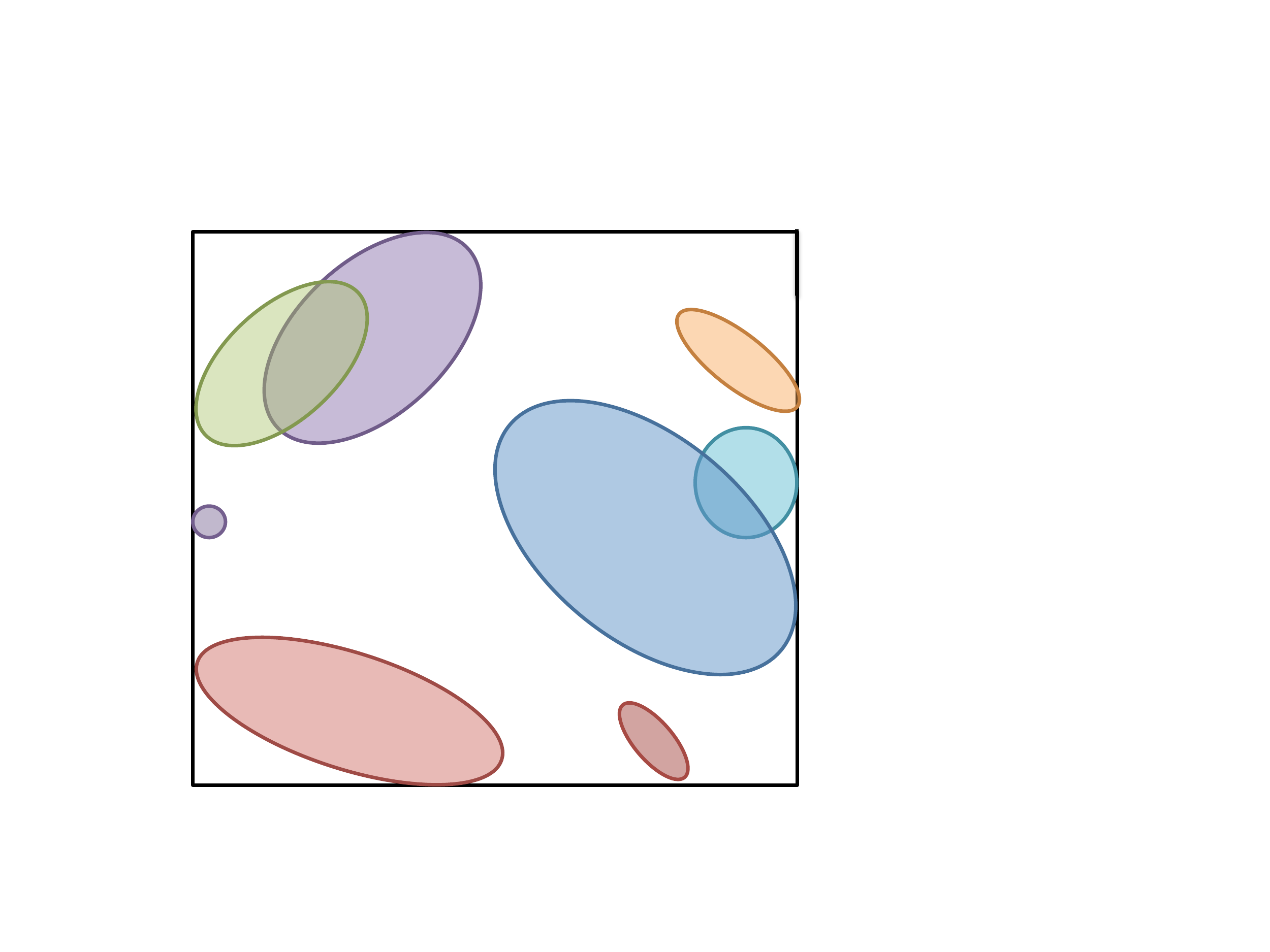}
\caption{No-slip Boundary}
\label{fig:1sub2}
\end{subfigure}
\caption{Boundary Conditions: (a) and (b) respectively shows treatment on periodic and no-slip boundaries. In (a), circles in the plot are supports of different modes. Mode 1 is separated into $A$, $B$, $C$, $D$ parts across boundary. $A$ stays at its location and $B$, $C$, $D$ are shifted to opposite sides($B'$, $C'$, $D'$). Mode 2 is separated into $E$, $F$ parts across boundary. $E$ stays at its location and $F'$ is shifted to opposite sides($F'$). In (b), supports of modes are restricted in the domain so that velocities are exactly $0$ on boundary. }\label{fig:digit}
\end{figure}       
For large enough flow domain, boundary condition for such turbulence generation has no influence on inner regions away from boundary. However, in order to retain flow properties near boundary, modes near boundary need to be treated differently. For periodic boundaries(Fig.\ref{fig:1sub1}), modes on boundary are separated into different parts and added to opposite boundaries.  Such treatment maintains exact same velocity value on opposite boundaries. 

For no-slip boundaries, modes are restricted to inner side of boundary(Fig.\ref{fig:1sub2}). In this way velocity and second-order moments on the boundary are exact zero. Also, characteristic length of modes near solid boundary is strictly restricted by its distance from boundary, which automatically creates damping effect near boundary.

\subsection{Input Spectrum}
\subsubsection{Isotropic case}
For high Reynolds number turbulence, homogeneous isotropic hypothesis is assumed to hold locally. Thus von Karman-Pao spectrum can be used to obtain spectral information of turbulence field. The von Karman-Pao spectrum is given by (\citet{juves1999stochastic}, \citet{saad2016scalable}):
       \begin{equation}\label{eq:22}
       E(k) = \alpha \frac{u'^2}{k_e}\frac{(k/k_e)^4}{[1+(k/k_e)^2]^{17/6}}\exp\Big[-2\Big(\frac{k}{k_{\eta}}\Big)^2\Big] 
       \end{equation}
       where $k_{\eta} = \epsilon^{1/4} \nu^{-3/4}$ is Kolmogorov length scale corresponding to viscous dissipation lengthscale. $\epsilon$ is turbulence dissipation rate from RANS data. $\alpha$ is determined from normalization of Eq.\ref{eq:10}: 
       \begin{equation}\label{eq:23}
       \alpha = \frac{55}{9\sqrt{\pi}}\frac{\Gamma(\frac{5}{6})}{\Gamma(\frac{1}{3})} \approx 1.453
       \end{equation}
       $k_e$ is the wavenumber related to most energetic eddies, could be determined by:
       \begin{equation}\label{eq:24}
       k_e = \sqrt{\pi}\frac{\Gamma(\frac{5}{6})}{\Gamma(\frac{1}{3})} \frac{1}{L} \approx \frac{0.746834}{L}
       \end{equation}
       $L = u'^3/\epsilon$ is integral lengthscale which could be computed from RANS data. 
\subsubsection{Wall turbulence case}
Eq.\ref{eq:20} defines a wavenumber $k_e$related to the energy containing structures, which could lead to a length scale $l_e$ related to $k_e$:
$$
l_e = \frac{2\pi}{k_e}
$$
$l_e$ corresponds to the size of most energetic eddies. In regions near wall, $l_e$ should not be larger that double the distance to the wall(\citet{shur2014synthetic}):
$$
l_e  \leqslant 2d_w
$$
where $d_w$ is the distance to wall. In regions far away from wall where damping effect is not important, expression of $k_e$ returns to isotropic von-Karman spectrum. Thus a modified expression of $k_e$ considering wall effect is as follow:
$$
k_e = \max(\sqrt{\pi}\frac{\Gamma(\frac{5}{6})}{\Gamma(\frac{1}{3})} \frac{1}{L} ,\frac{\pi}{d_w})
$$
\subsection{Spatial-Spectral decomposition}
\begin{figure}[htbp]
\centering
\begin{subfigure}{.4\textwidth}
\centering
\includegraphics[width=1.1\linewidth]{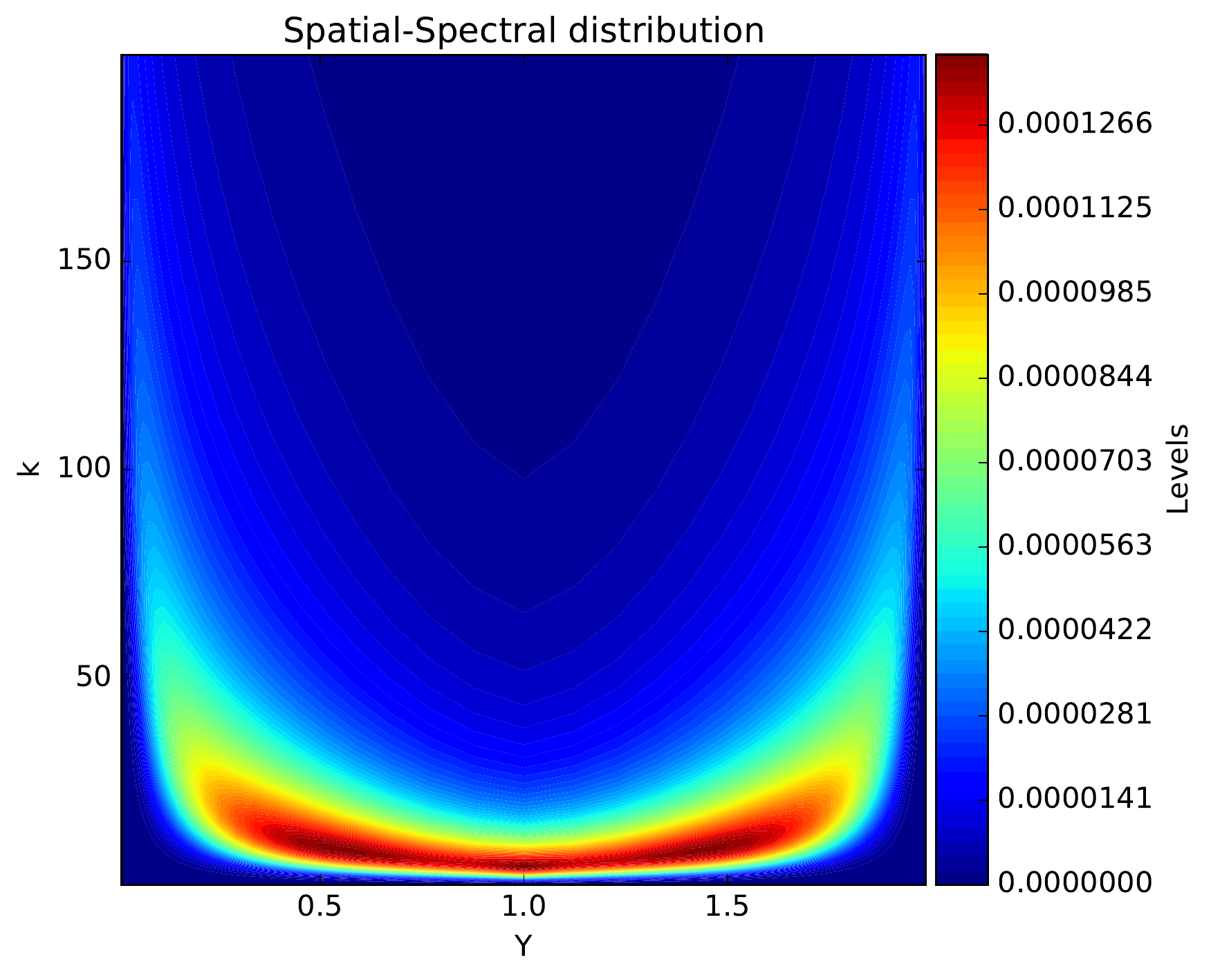}
\caption{spatial-spectral distribution of TKE}
\label{fig:3sub1}
\end{subfigure}%
\begin{subfigure}{.4\textwidth}
\centering
\includegraphics[width=0.9\linewidth]{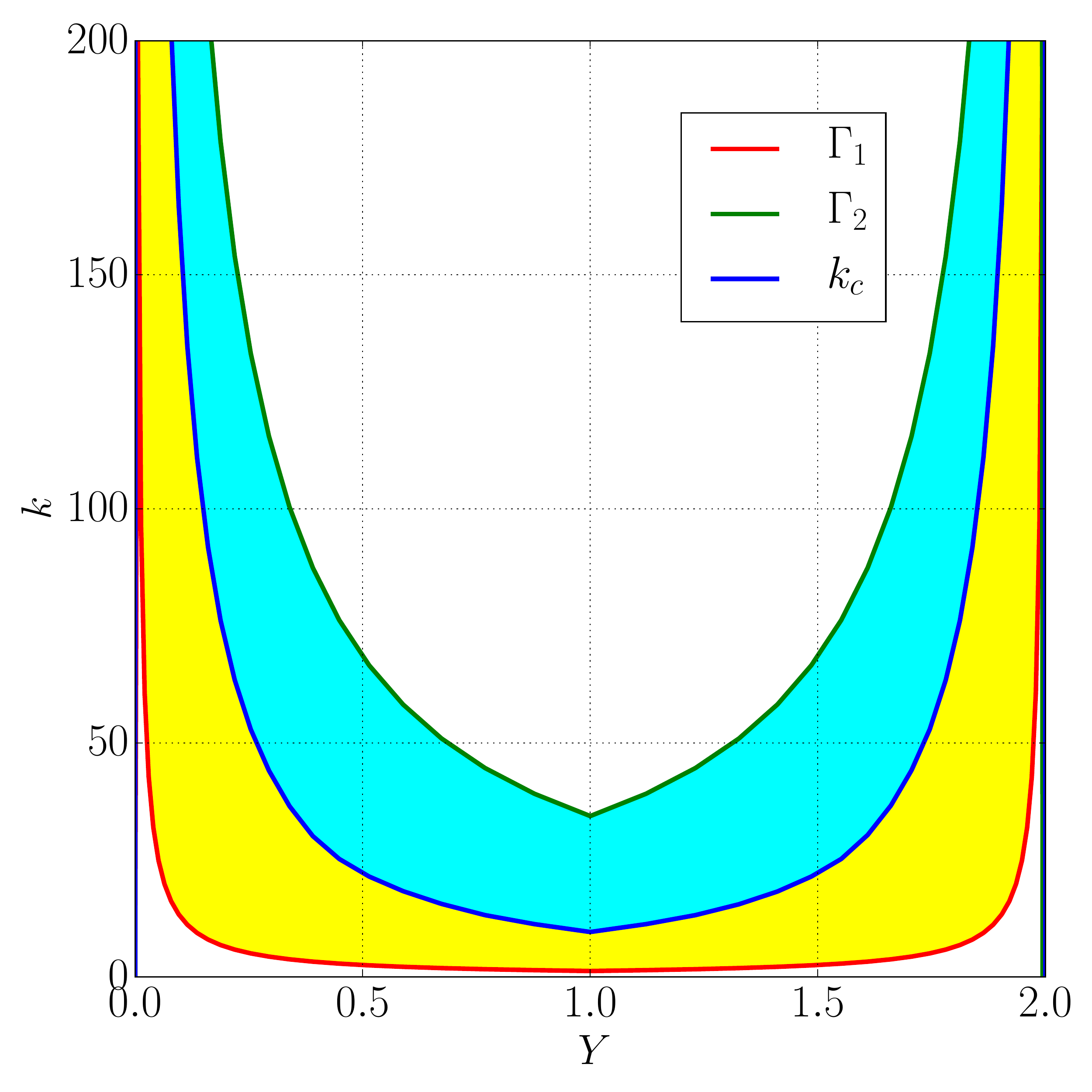}
\caption{Cut off boundaries}
\label{fig:3sub2}
\end{subfigure}
\caption{Channel flow: (a) presents TKE distribution in Physical-Spectral space. (b) shows the regions and cutoff bounds for turbulence generation. Generation process is only conducted in yellow and blue regions using different time-advance procedures proposed in Sec.\ref{Generation of dynamic turbulence field}. Yellow part is large scale generation with smaller wavenumber. Blue part is small scale with large wavenumber. White-colored area is discarded in generation process. }
\label{fig:3}
\end{figure}

Eq.\ref{eq:3} gives a decomposition of velocity potential field in both physical and spectral space, which allows a large reduction of computational cost of reconstruction of turbulence field, especially for anisotropic inhomogeneous turbulence. Consider a fully developed channel flow the simulation result of which is shown in Section  3.2(Fig.\ref{fig:2}). 
Such flow is basically 1 dimensional, in which all turbulence quantities are only functions of $y$. Eq.\ref{eq:19} defines the turbulent kinetic energy(TKE) distribution in spectral space, i.e. $E(k,k_t)$. Thus, a spatial-spectral distribution of TKE is defined as follow:
$$
E(k,y) = E(k,k_t(y))
$$
which satisfies the following normalization condition:
$$
k_t(y) = \int_0^{\infty} E(k,k_t(y))\dif k
$$
Such distribution of channel flow (Fig.\ref{fig:3}) gives a special insight of the energy distribution and flow structure of turbulence field. In most part of the flow field which is far away from wall, most of turbulent kinetic energy concentrate in a very narrow area of spectral space, which only consists a small portion of computational cost in the simulation. In turbulent boundary layer near the wall, the distribution of TKE in spectral space becomes very wide and considerably increases computation cost. Based on the generation method in Section 2.1, a scale-reduction algorithm is designed to cut off unnecessary computation while still capture the energetic structures in the flow field. Define error $e$ representing the energy lost ratio in order to reduce computational cost:
$$
e = 1-\frac{k_t^{e}}{k_t}
$$
Define two boundary curves $\Gamma_1(y)$, $\Gamma_2(y)$ as follow:
$$
\int_{\Gamma_1^e(y)}^{\Gamma_2^e(y)}E(k,y)\dif k = k_t^e(y)
$$
The generation process only need to be operated in a small band of physical-spectral space (Fig.\ref{fig:3}) between $\Gamma_1$ and $\Gamma_2$ to reconstruct most of TKE up to an error $e$.

\section{Numerical Results}\label{Numerical Results}
\subsection{Isotropic homogeneous turbulence}
\subsubsection{Spatial structure}
\begin{figure}[htbp]
\centering
\includegraphics[width=0.5\linewidth]{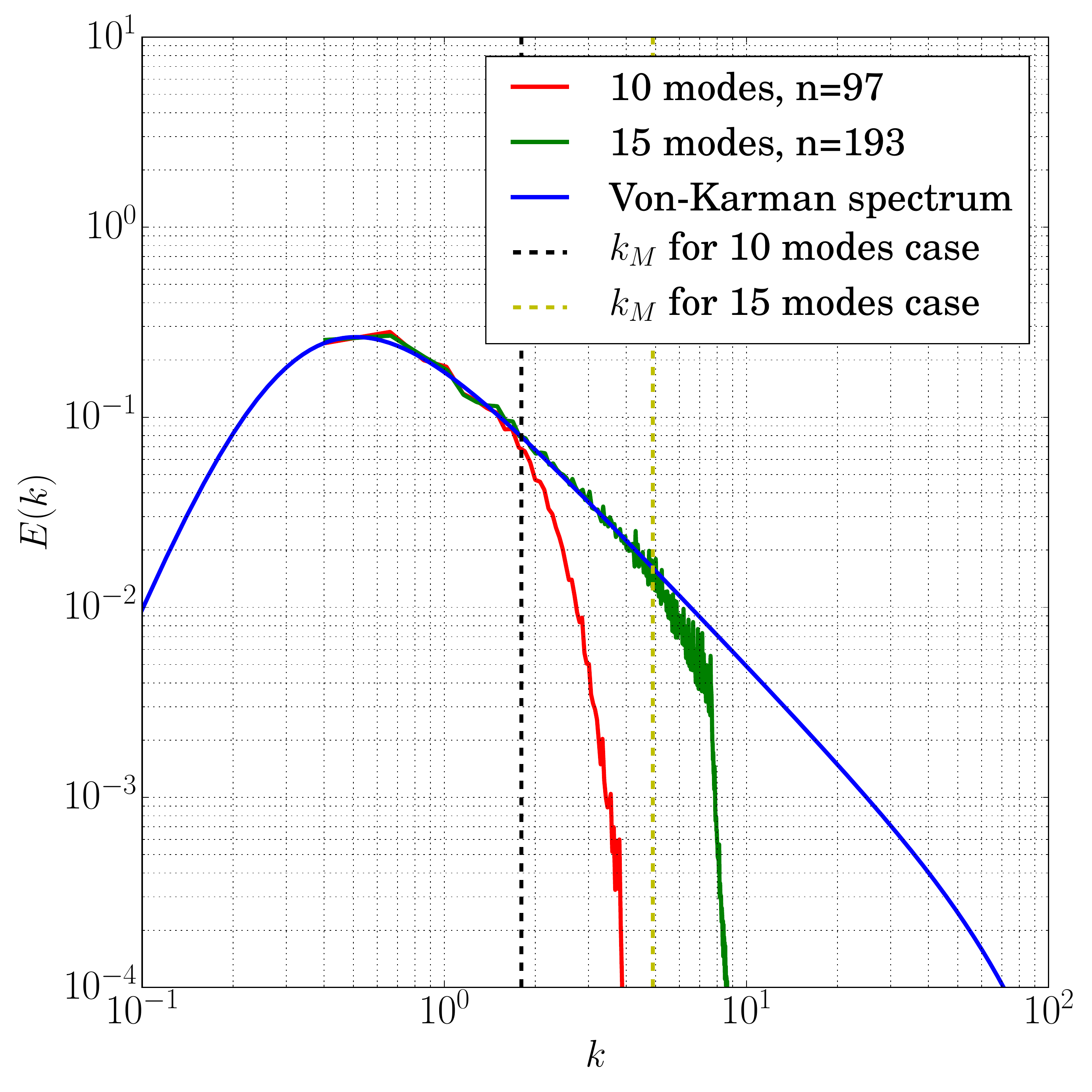}
\caption{Spectrum of generated turbulence: In this simulation $k_t=0.5$, $\nu=7\times 10^{-4}$, $\epsilon=0.0849$ are nondimensional parameters from previous DNS results. $k_{min}$ is the cutoff wavenumber of original simulation representing the largest scale. $k_0 = 0.3$, $d_k=0.2$ for wavenumber generation. The result shows good agreement with input spectrum within the range of wavenumbers of modes used in the generation. With mode number increased, the spectrum range that could be accurately captured gets larger. $k_M$ is the largest wavenumber used in generation. Normally $k_M\leqslant k_N$. $k_N$ is Nyquist wavenumber of mesh. }\label{fig:4}
\end{figure}

\begin{figure}[htbp]
\centering
\begin{subfigure}{.5\textwidth}
\centering
\includegraphics[width=1.0\linewidth]{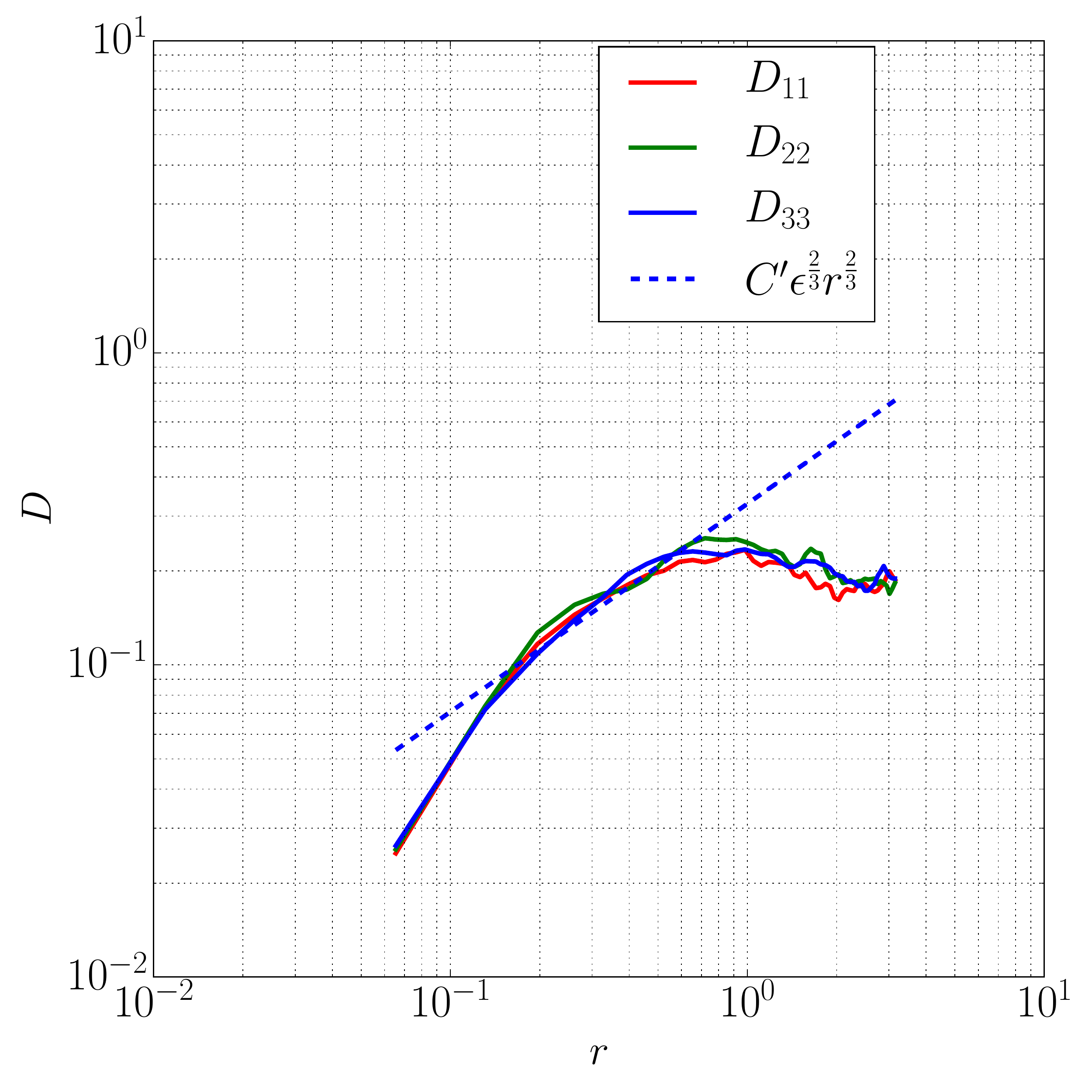}
\caption{10modes, n=97}
\label{fig:5sub1}
\end{subfigure}%
\begin{subfigure}{.5\textwidth}
\centering
\includegraphics[width=1.0\linewidth]{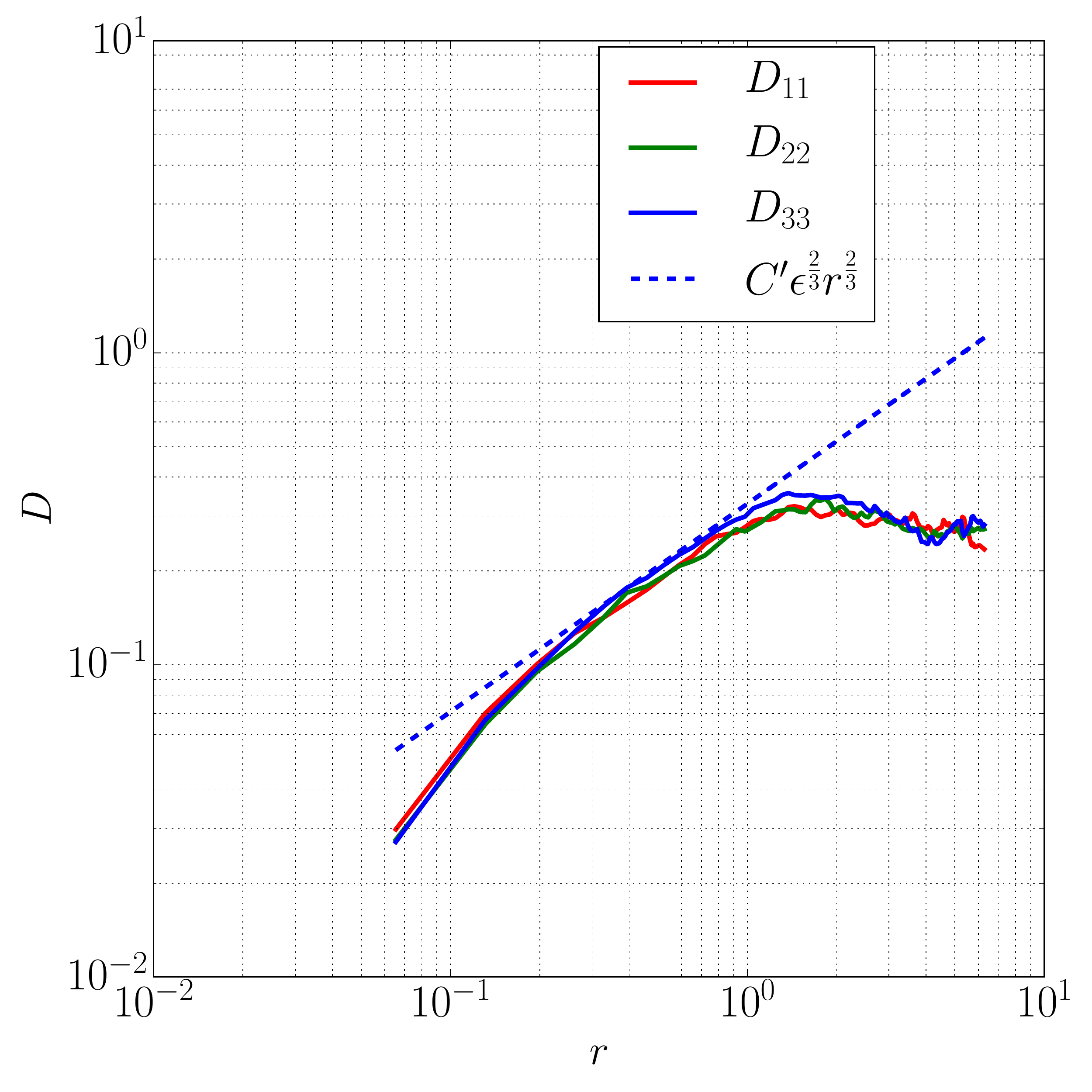}
\caption{15 modes, n=193}
\label{fig:5sub2}
\end{subfigure}
\caption{Second order structure functions of simulations with different number of modes and grid points. Dashed lines are theoretical results with constant $C'=1.7$. The simulation shows good agreement compared with theoretical results in inertial subrange }
\label{fig:5}
\end{figure}

\begin{figure}[htbp]
\centering
\begin{subfigure}{.5\textwidth}
\centering
\includegraphics[width=1.0\linewidth]{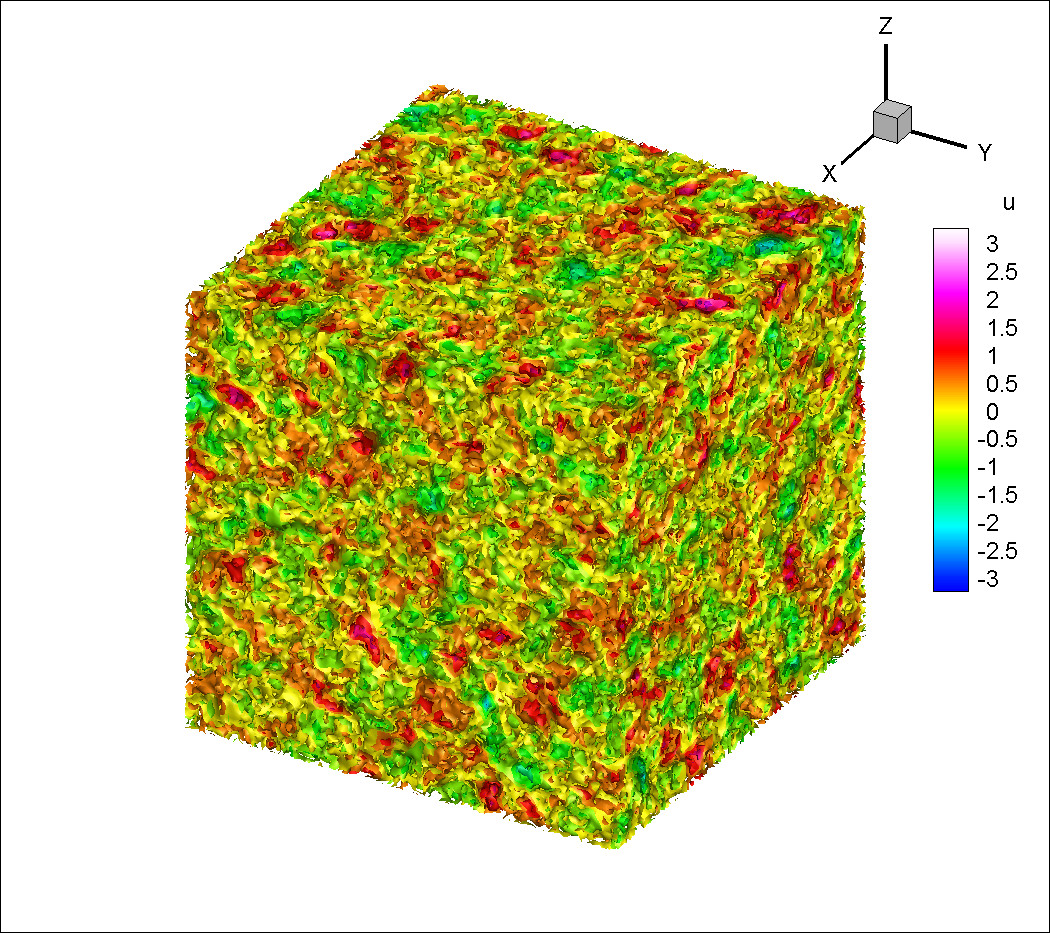}
\caption{Isosurface of $u$}
\label{fig:6sub1}
\end{subfigure}%
\begin{subfigure}{.5\textwidth}
\centering
\includegraphics[width=1.0\linewidth]{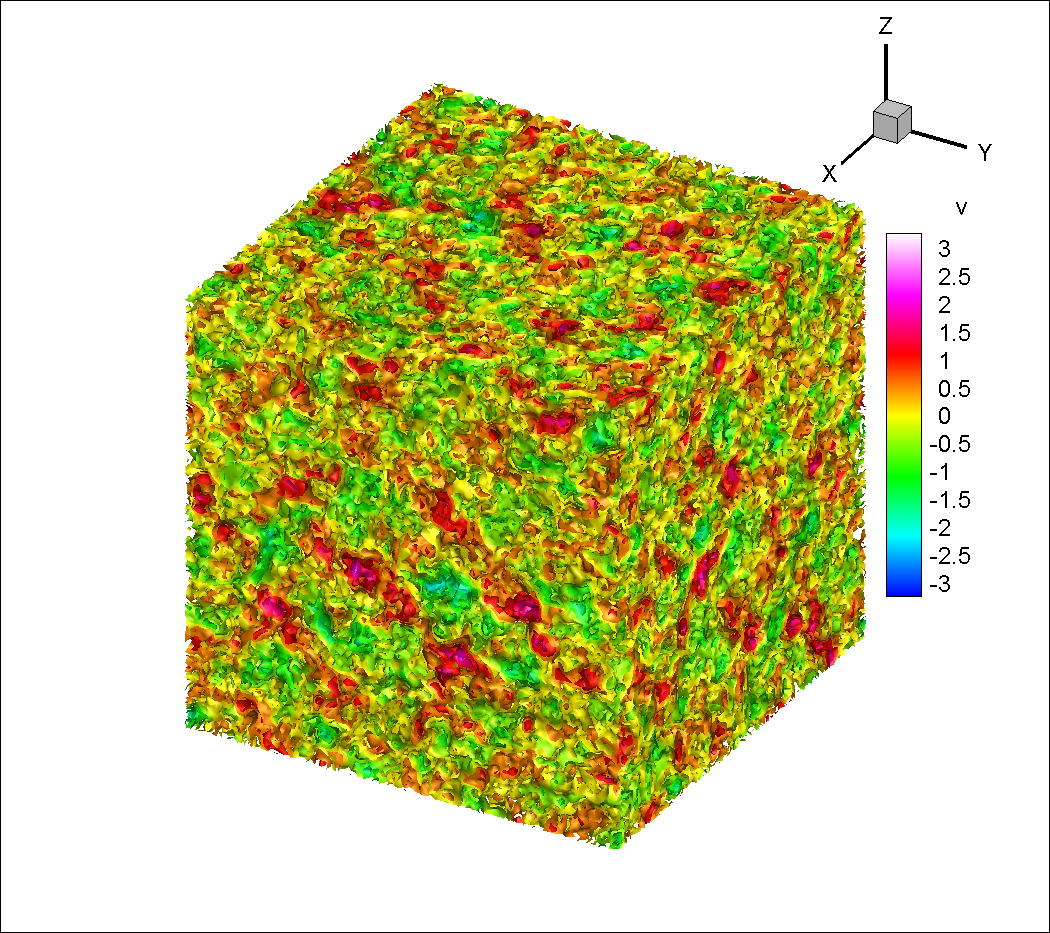}
\caption{Isosurface of $v$}
\label{fig:6sub2}
\end{subfigure}
\begin{subfigure}{.5\textwidth}
\centering
\includegraphics[width=1.0\linewidth]{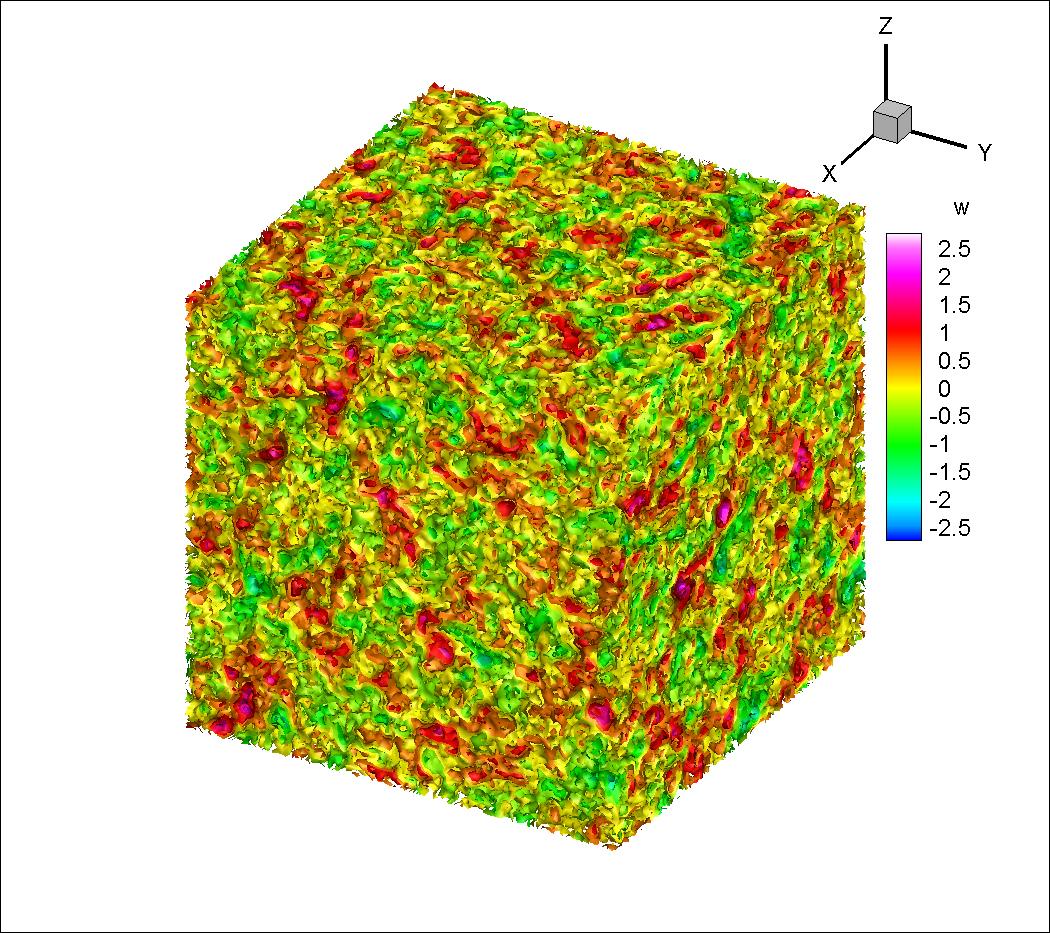}
\caption{Isosurface of $w$}
\label{fig:6sub3}
\end{subfigure}%
\begin{subfigure}{.5\textwidth}
\centering
\includegraphics[width=1.0\linewidth]{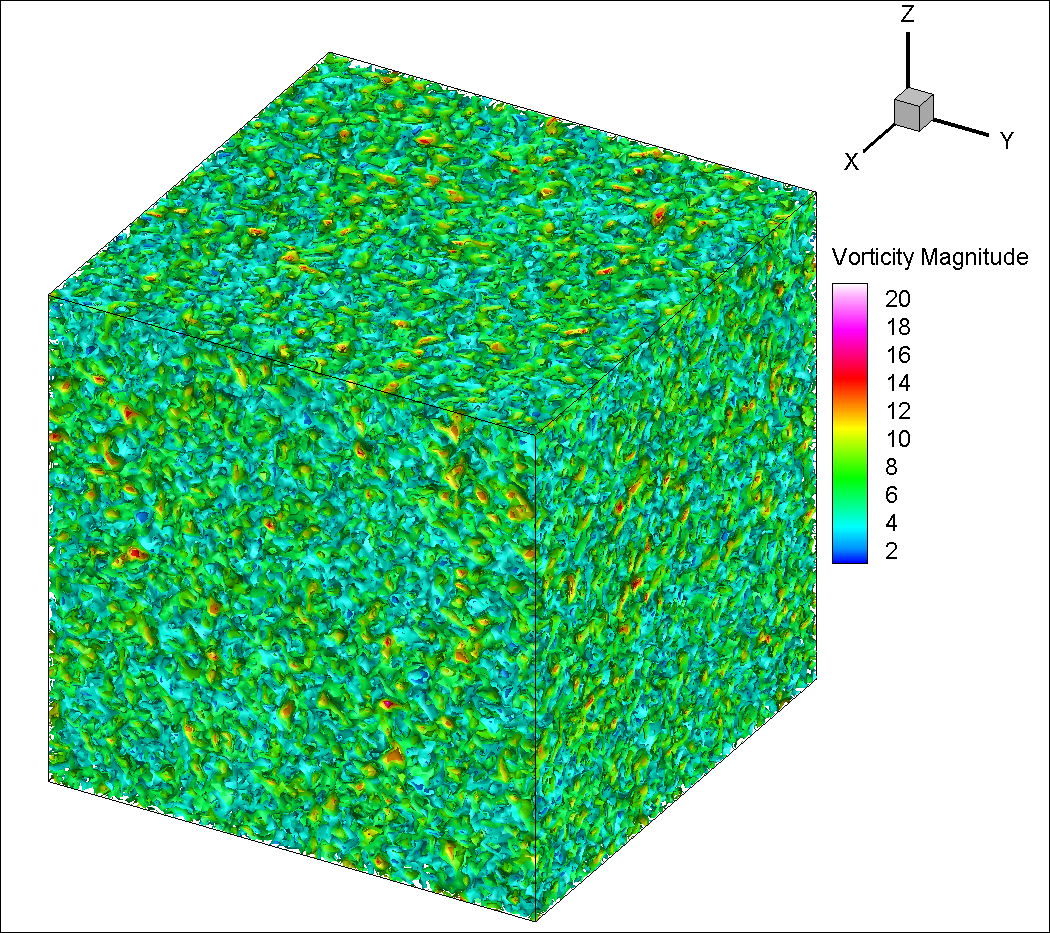}
\caption{Isosurface of vorticity magnitude}
\label{fig:6sub4}
\end{subfigure}
\caption{Isosurface Results of Isotropic Homogeneous Turbulence Generation}\label
{fig:6}
\end{figure}

Generation of isotropic homogeneous turbulence is an important way to validate various properties of turbulence synthesis models. An isotropic homogeneous turbulence case is computed in order to verify the model constructed in Methodology section. Because RANS type model cannot compute $k_t$ and $\epsilon$ of isotropic homogeneous turbulence, such data  is obtained from previous DNS results (\citet{kaneda2003energy}). Spectrum of generated turbulence with different modes and different resolution are compared with von-Karman spectrum(Fig.\ref{fig:4}). Wavenumbers of modes are generated with the following formula(\citet{juves1999stochastic}):
$$
k_n = \exp[\ln k_0 +nd_k],  = 0,1,2,.....M
$$ 
Where $k_0$ is the first wavenumber of the sequence, $d_k$ is a parameter to control the distances between wavenumbers. 
Another quantity that could be used to exam spatial structure of generated turbulence is structure function defined as follow:
$$
D_{11} (r,0,0) = \langle[u(x+r,y,z,t)-u(x,y,z,t)]^{2}\rangle
$$
$$
D_{22} (0,r,0) = \langle[u(x,y+r,z,t)-u(x,y,z,t)]^{2}\rangle
$$
$$
D_{33} (0,0,r) = \langle[u(x,y,z+r,t)-u(x,y,z,t)]^{2}\rangle
$$
From previous theoretical and experimental research(\citet{fung1992kinematic},\citet{ishihara2009study},etc.), second order structure function has the following form in inertial subrange:
$$
D_{11}=D_{22}=D_{33}=C'\epsilon^{\frac{2}{3}}r^{\frac{2}{3}}
$$
where C' is a constant. In \citet{fung1992kinematic} the value of $C'$ is equal to 1.7. Numerical results of $D_{11}$, $D_{22}$ and $D_{33}$ are compared  with theoretical solutions(Fig.\ref{fig:5}). The energy spectrum results and second order structural function results indicates that this turbulence generation method gives right spatial turbulence structure in homogeneous isotropic cases. The iso-surface of the numerical results are shown in Fig.\ref{fig:6}. 

\subsubsection{Time correlation}
\begin{figure}[htbp]
\centering
\begin{subfigure}{.5\textwidth}
\centering
\includegraphics[width=1.0\linewidth]{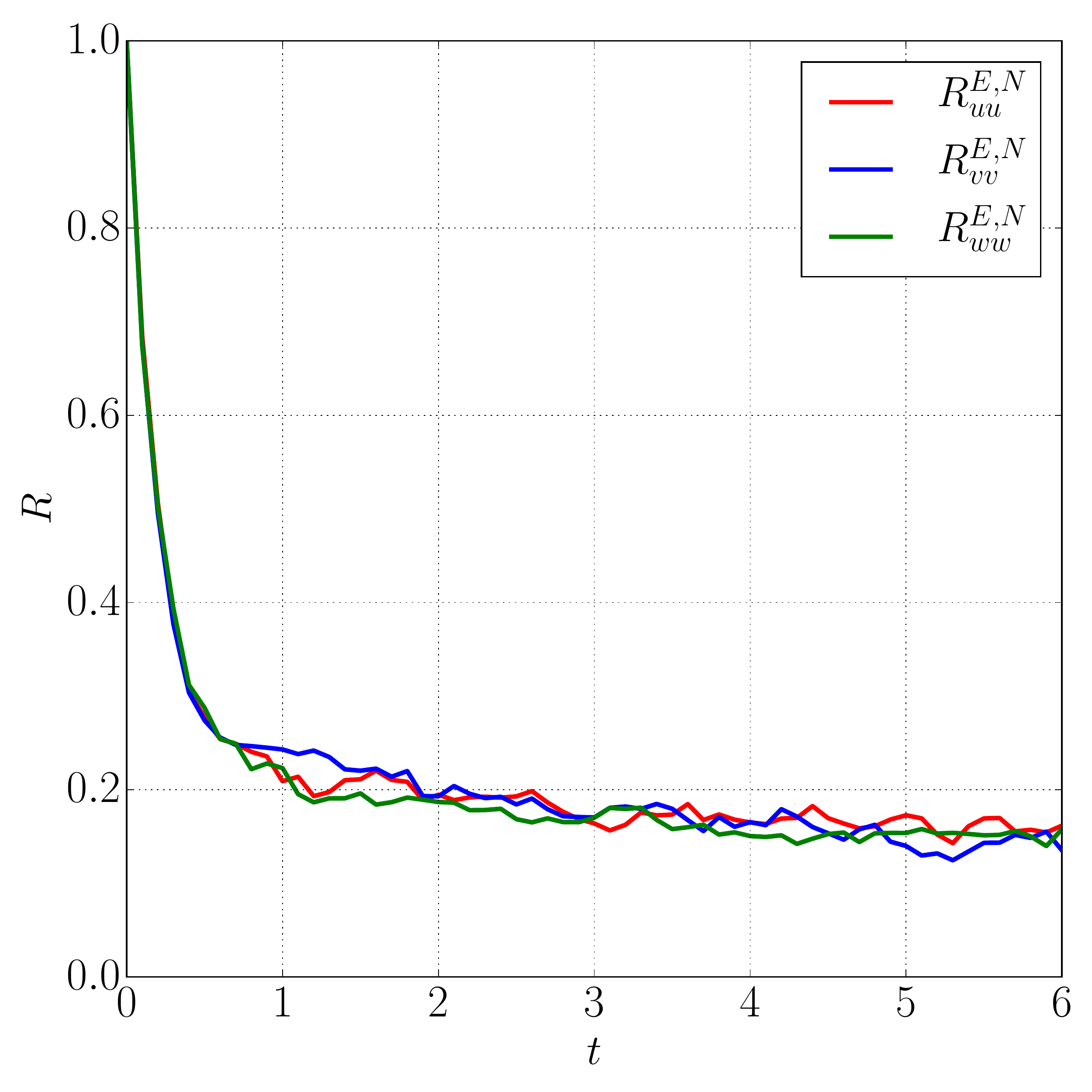}
\caption{Normalized Eulerian autocorrelation}
\label{fig:7sub1}
\end{subfigure}%
\begin{subfigure}{.5\textwidth}
\centering
\includegraphics[width=1.0\linewidth]{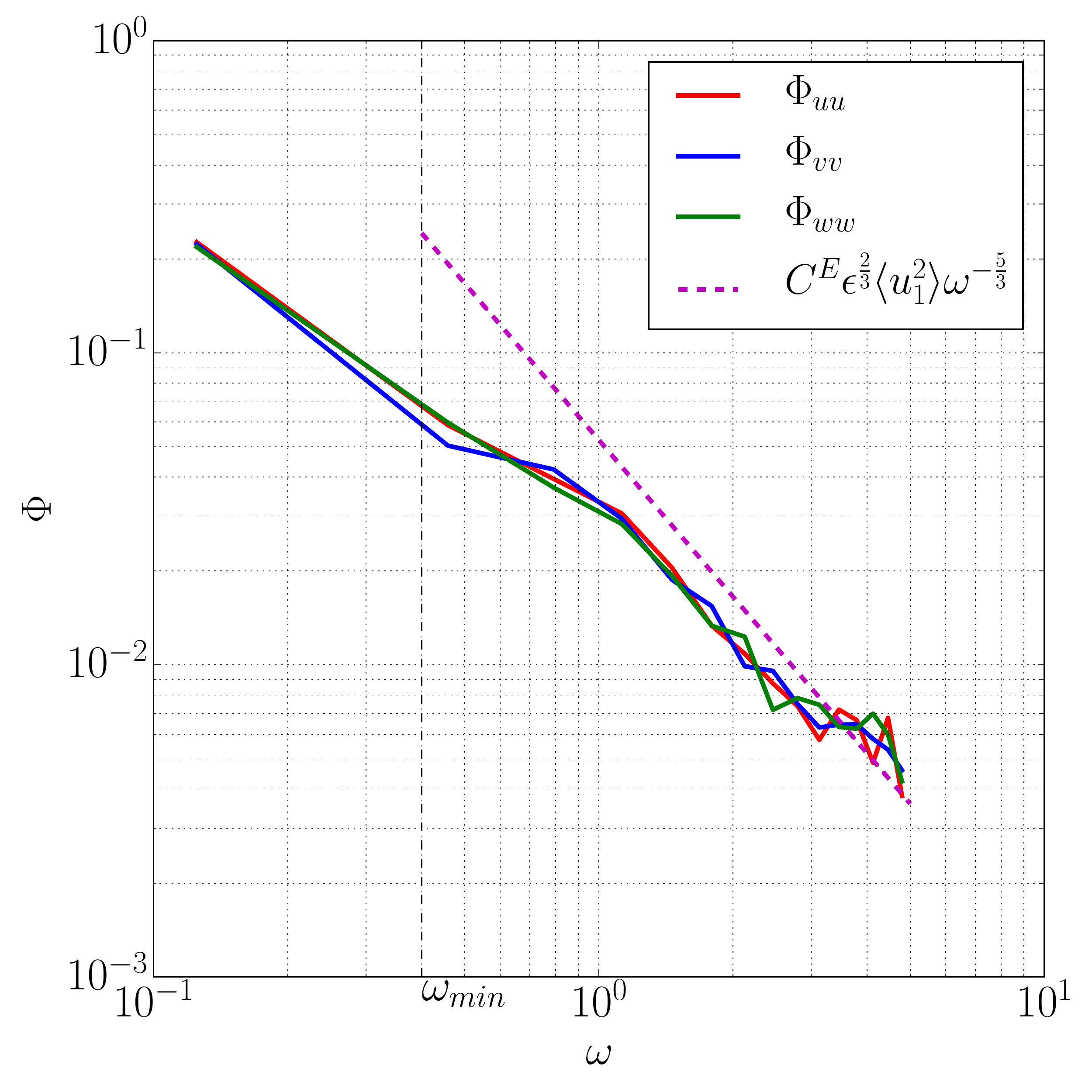}
\caption{Eulerian frequency spectrum}
\label{fig:7sub2}
\end{subfigure}
\caption{Eulerian temporal and frequency properties of simulation results. Frequency spectrums are compared with analytical results from previous study. $\omega_{min}$ is the lower boundary of frequency inertial subrange. It can be observed that the Eulerian frequency spectrums of numerical results approximate theoretical solution in inertial subrange, especially at high frequency range. }\label{fig:7}
\end{figure}

Eulerian autocorrelation is defined as follow:
$$
R^{E}_{uu}(\tau) = \langle u(t)u(t+\tau)\rangle
$$
$$
R^{E}_{vv}(\tau) = \langle v(t)v(t+\tau)\rangle
$$
$$
R^{E}_{ww}(\tau) = \langle w(t)w(t+\tau)\rangle
$$
Normalized Eulerian autocorrelation can be computed as follow:
$$
R^{E,N}_{11}(\tau) = \frac{R^{E}_{uu}(\tau)}{R^{E}_{uu}(0)}
$$
$$
R^{E,N}_{22}(\tau) = \frac{R^{E}_{vv}(\tau)}{R^{E}_{vv}(0)}
$$
$$
R^{E,N}_{33}(\tau) = \frac{R^{E}_{ww}(\tau)}{R^{E}_{ww}(0)}
$$
Eulerian frequency spectrums are defined as the Fourier transform of Eulerian autocorrelation:
$$
\Phi^{E}_{uu}(\omega) = \int_{\mathbb{R}}R^{E}_{uu}(\tau) \eup ^{-i\omega \tau}\dif \tau 
$$
$$
\Phi^{E}_{vv}(\omega) = \int_{\mathbb{R}}R^{E}_{vv}(\tau) \eup ^{-i\omega \tau}\dif \tau 
$$
$$
\Phi^{E}_{ww}(\omega) = \int_{\mathbb{R}}R^{E}_{ww}(\tau) \eup ^{-i\omega \tau}\dif \tau 
$$
\citet{fung1992kinematic},\citet{ishihara2009study} suggest that for isotropic homogeneous turbulence, Eulerian frequency spectrum in inertial subrange can be approximated as follow:
$$
\Phi^{E}_{uu}(\omega) = \Phi^{E}_{vv}(\omega) = \Phi^{E}_{ww}(\omega) \approx C^{E}\epsilon^{\frac{2}{3}}\langle u^{2}_{1}\rangle\omega^{-\frac{5}{3}}
$$
where $C^{E}=0.46$ is a constant from \citet{ishihara2009study}. Results of Eulerian autocorrelation and frequency spectrum are shown in Fig.\ref{fig:7}. 

\subsection{Fully Developed Turbulent Channel Flow}
\begin{figure}[htbp]
\centering
\begin{subfigure}{.5\textwidth}\label{fig:8sub1}
\centering
\includegraphics[width=1.0\linewidth]{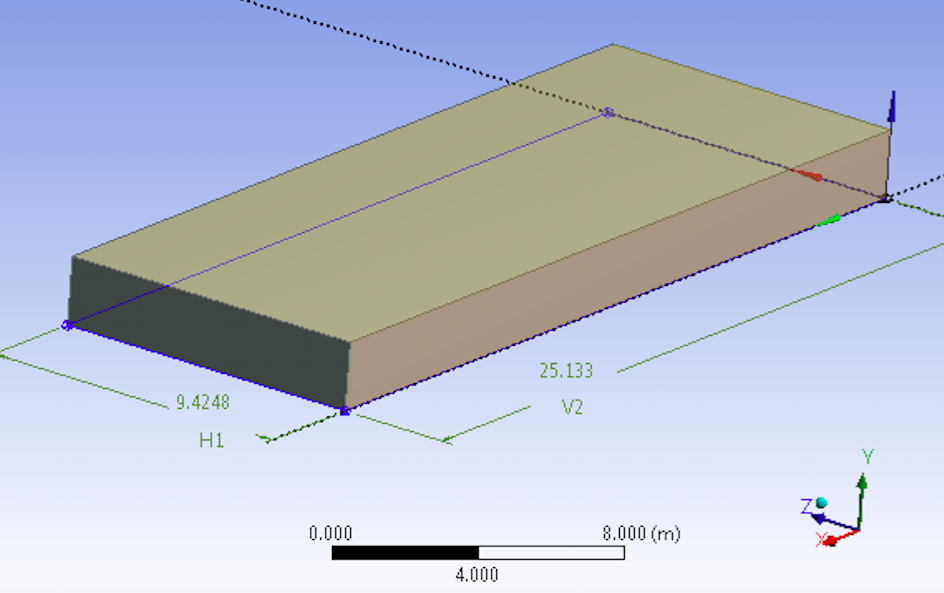}
\caption{Geometry setup}
\end{subfigure}%
\begin{subfigure}{.5\textwidth}\label{fig:8sub2}
\centering
\includegraphics[width=1.0\linewidth]{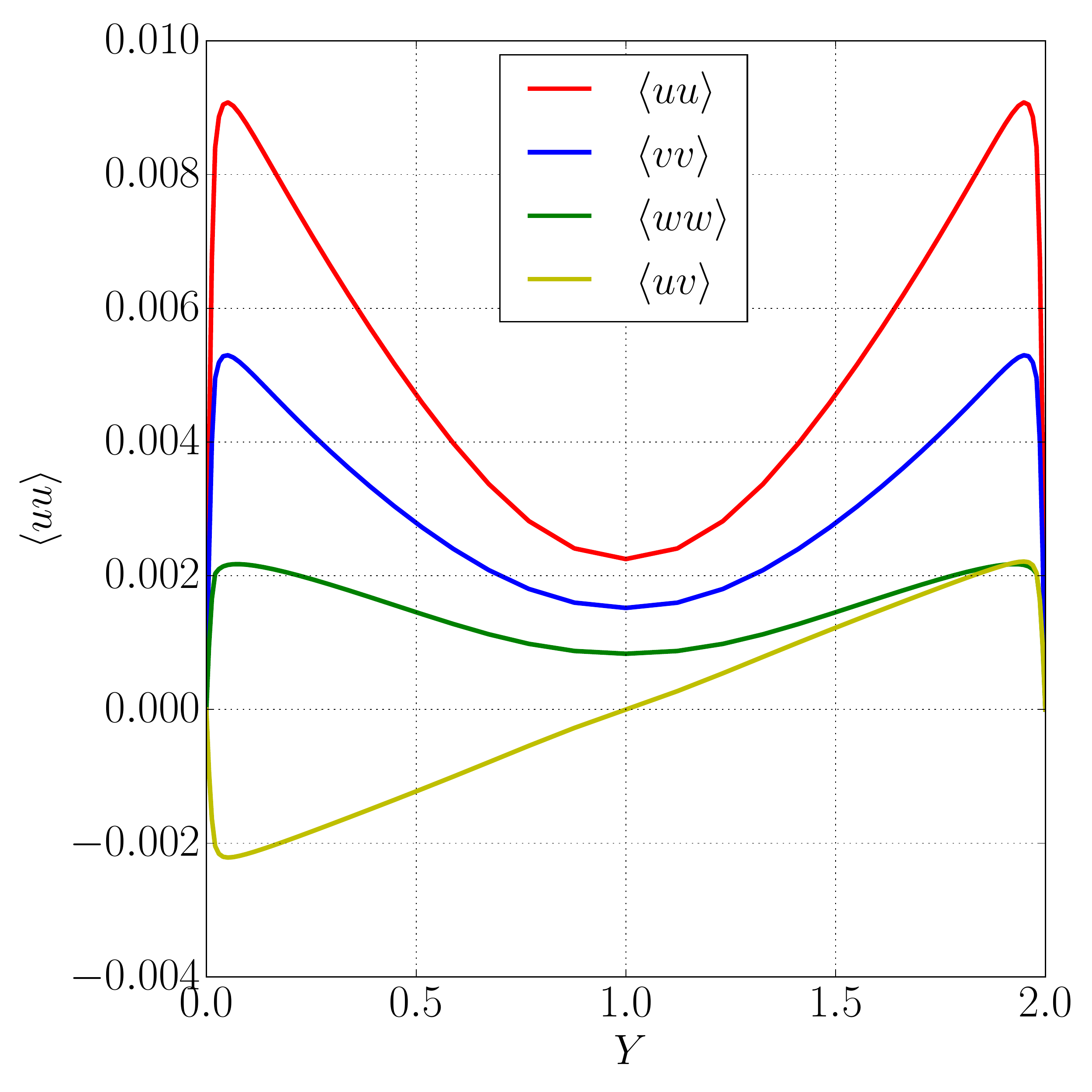}
\caption{RANS data}
\end{subfigure}
\caption{Channel flow: (a) is geometry setup of this channel flow, X is the streamwise direction of the channel flow. $Z=0$ and $Z=3\pi$ are periodic boundaries. $Y = 0$ and $Y=2$ are no-slip boundaries. Gradient of mean turbulence field only exist on $Y$ direction, while $X$ and $Z$ direction are uniform. (b) is RANS data of this channel flow from Reynolds Stress Model simulation. }\label{fig:8}
\end{figure}

RANS data of channel flow case in Fig.2 was computed using Reynolds Stress Model to obtain full Reynolds stress and turbulence dissipation rate(Fig.\ref{fig:8}). Turbulent boundary layer was fully resolved, including several grid points in viscous sublayer. This RANS data is used as input data for turbulence generation. 
Residual of RSM simulation results are in Table.\ref{table_Residual}. 
\begin{table}[!h]
\caption{Residuals}
\label{table_Residual}
\begin{subtable}{.5\textwidth}
\begin{tabular}{llllll}
\hline
Equation &Continuity&$U$ momentum&$V$ momentum&$W$ momentum&$k$ equation\\
\hline
Residual &$1.23\times 10^{-6}$&$2.04\times 10^{-6}$&$7.04\times 10^{-9}$&$6.61\times 10^{-9}$&$2.3044 \times 10^{-6}$ \\\hline
\end{tabular}
\vspace*{-4pt}
\end{subtable}

\begin{subtable}{.5\textwidth}
\begin{tabular}{llllll}
\hline
Equation & $\epsilon$ &$\langle uu\rangle$&$\langle vv\rangle$&$\langle ww\rangle$ &$\langle uv\rangle$\\
\hline
Residual &$5.81\times 10^{-6}$&$5.59\times 10^{-6}$&$5.72\times 10^{-6}$&$5.68\times 10^{-6}$&$5.44 \times 10^{-6}$ \\\hline
\end{tabular}
\vspace*{-4pt}
\end{subtable}
\end{table}
 
\begin{figure}[htbp]
\centering
\begin{subfigure}{.5\textwidth}
\centering
\includegraphics[width=0.7\linewidth]{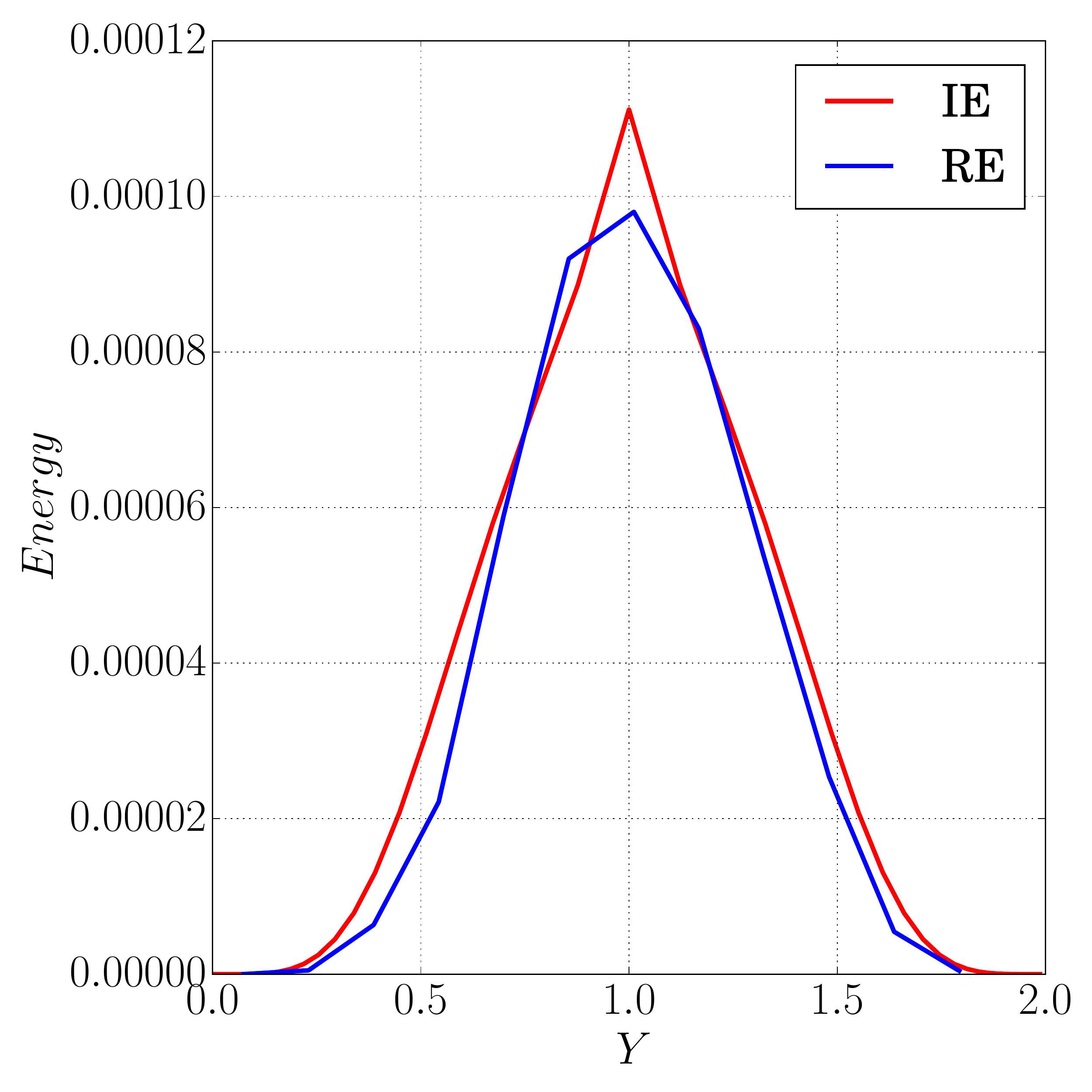}
\caption{$k=3.0$}
\label{fig:9sub1}
\end{subfigure}%
\begin{subfigure}{.5\textwidth}
\centering
\includegraphics[width=0.7\linewidth]{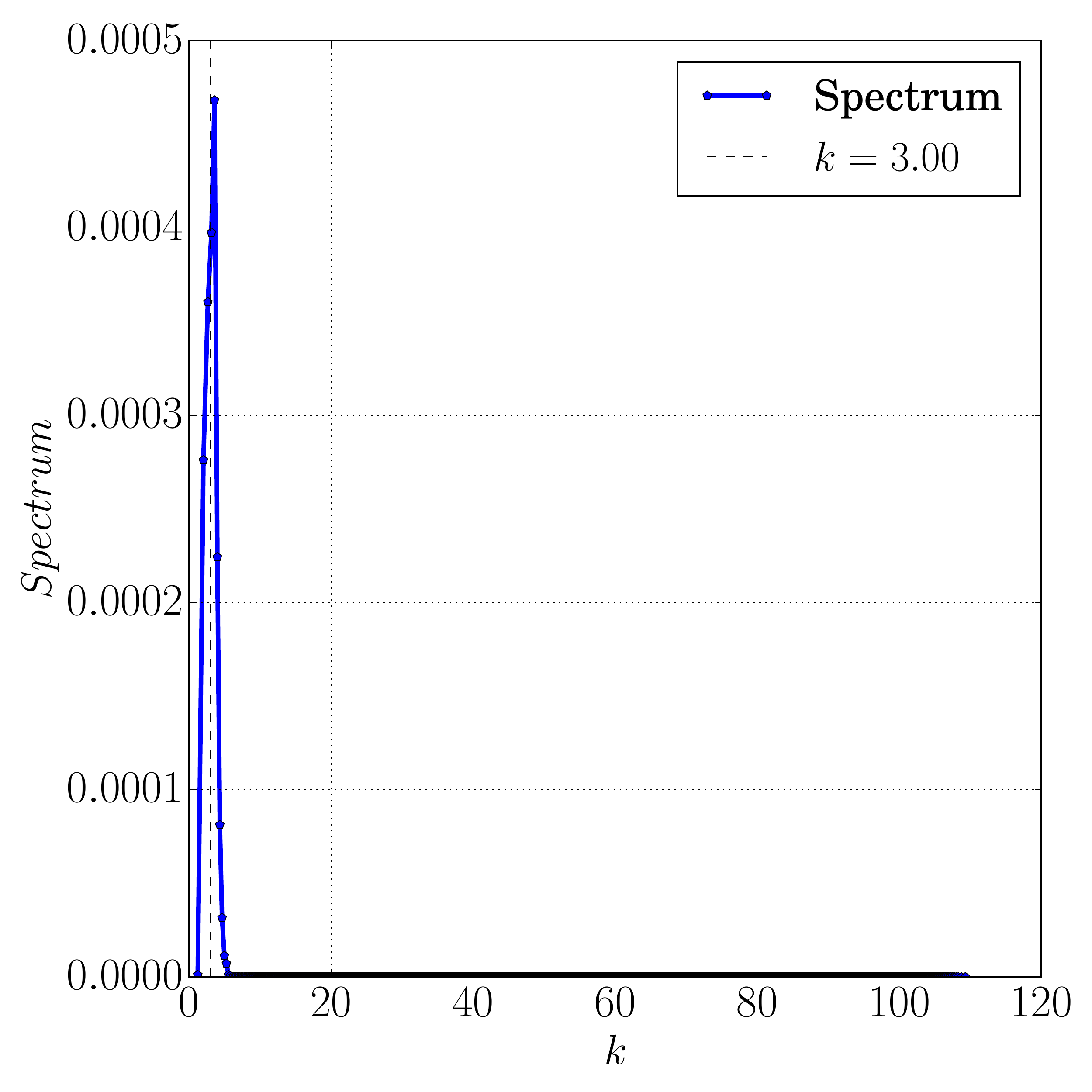}
\caption{$k=3.0$}
\label{fig:9sub2}
\end{subfigure}

\begin{subfigure}{.5\textwidth}
\centering
\includegraphics[width=0.7\linewidth]{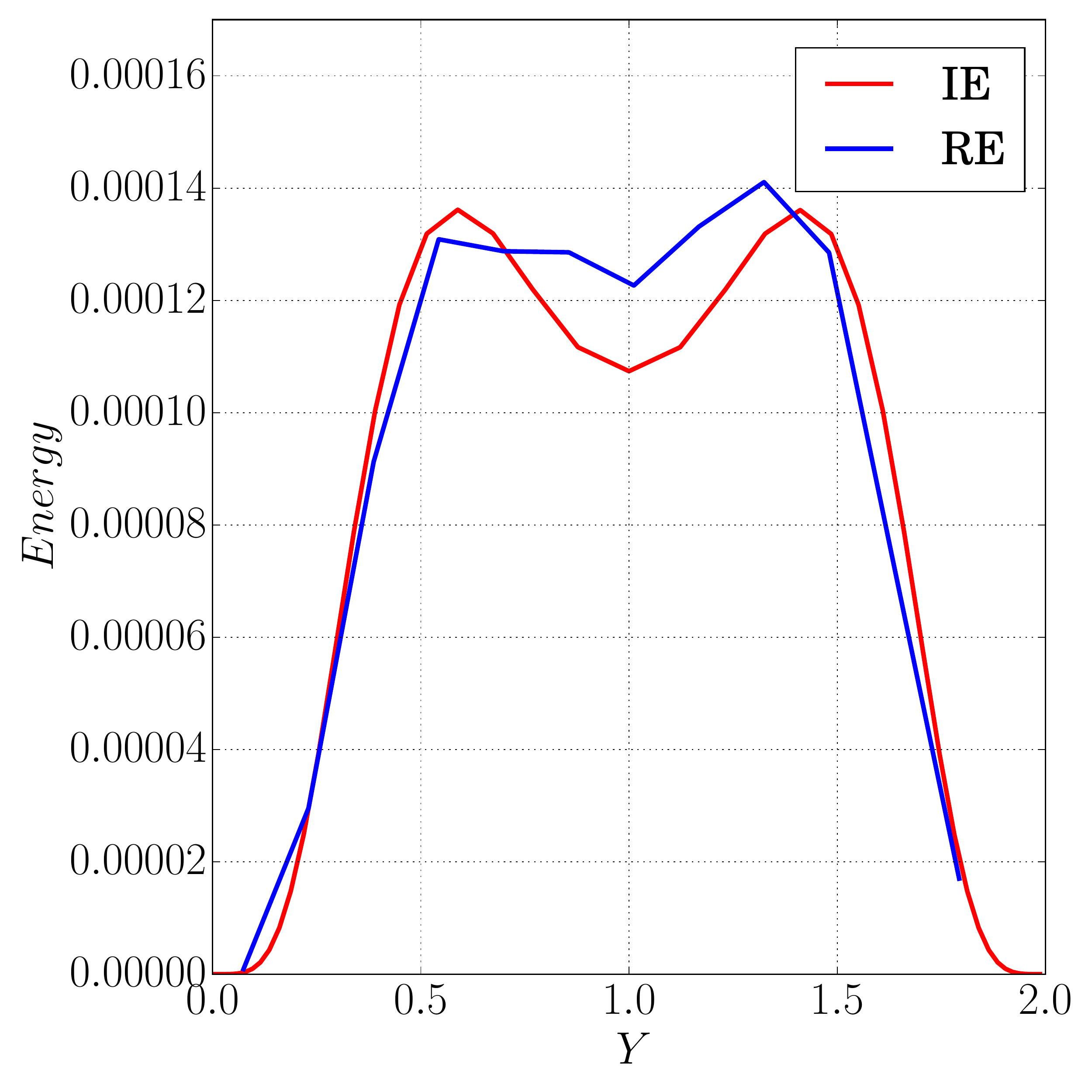}
\caption{$k=7.38$}
\label{fig:3sub3}
\end{subfigure}%
\begin{subfigure}{.5\textwidth}
\centering
\includegraphics[width=0.7\linewidth]{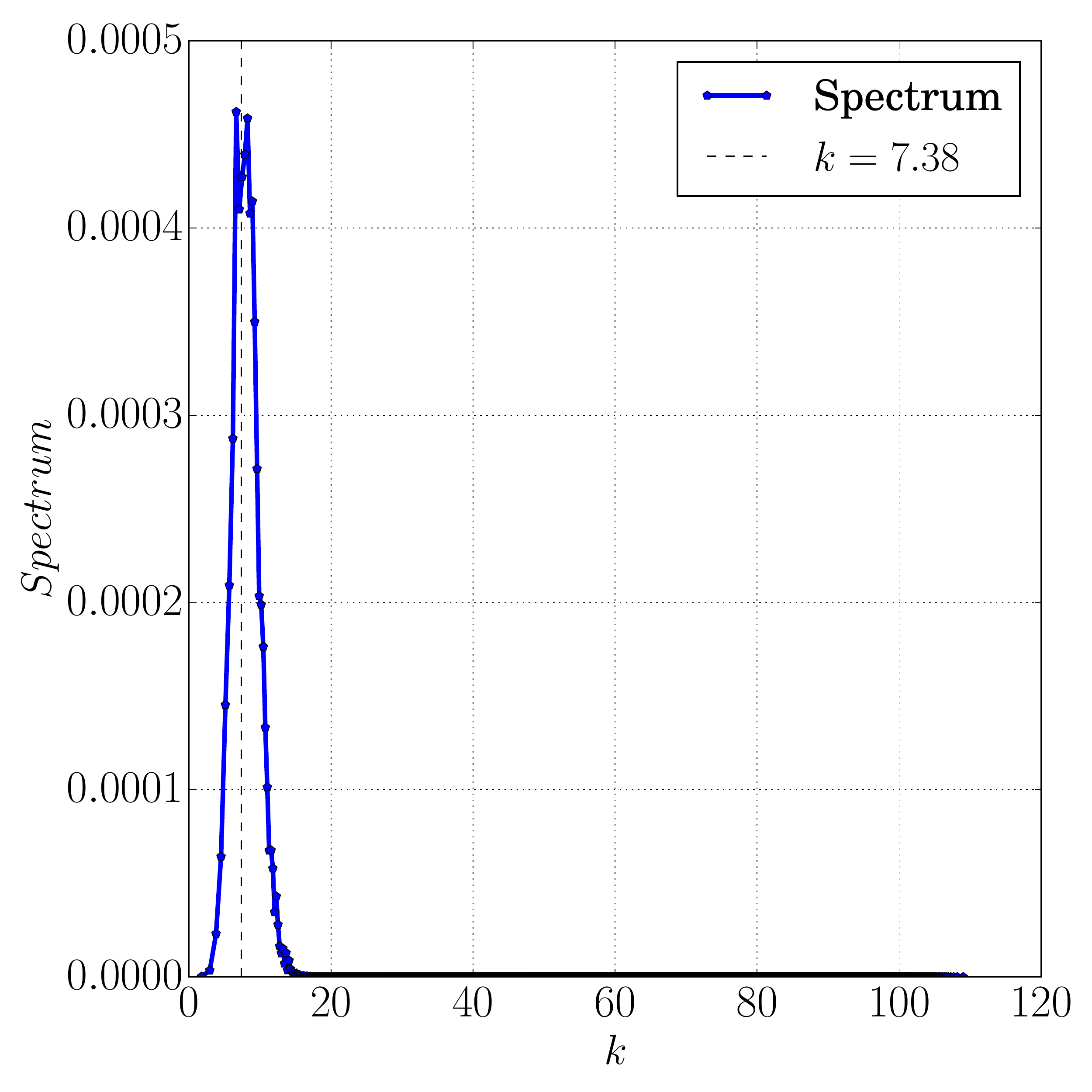}
\caption{$k=7.38$}
\label{fig:9sub4}
\end{subfigure}

\begin{subfigure}{.5\textwidth}
\centering
\includegraphics[width=0.7\linewidth]{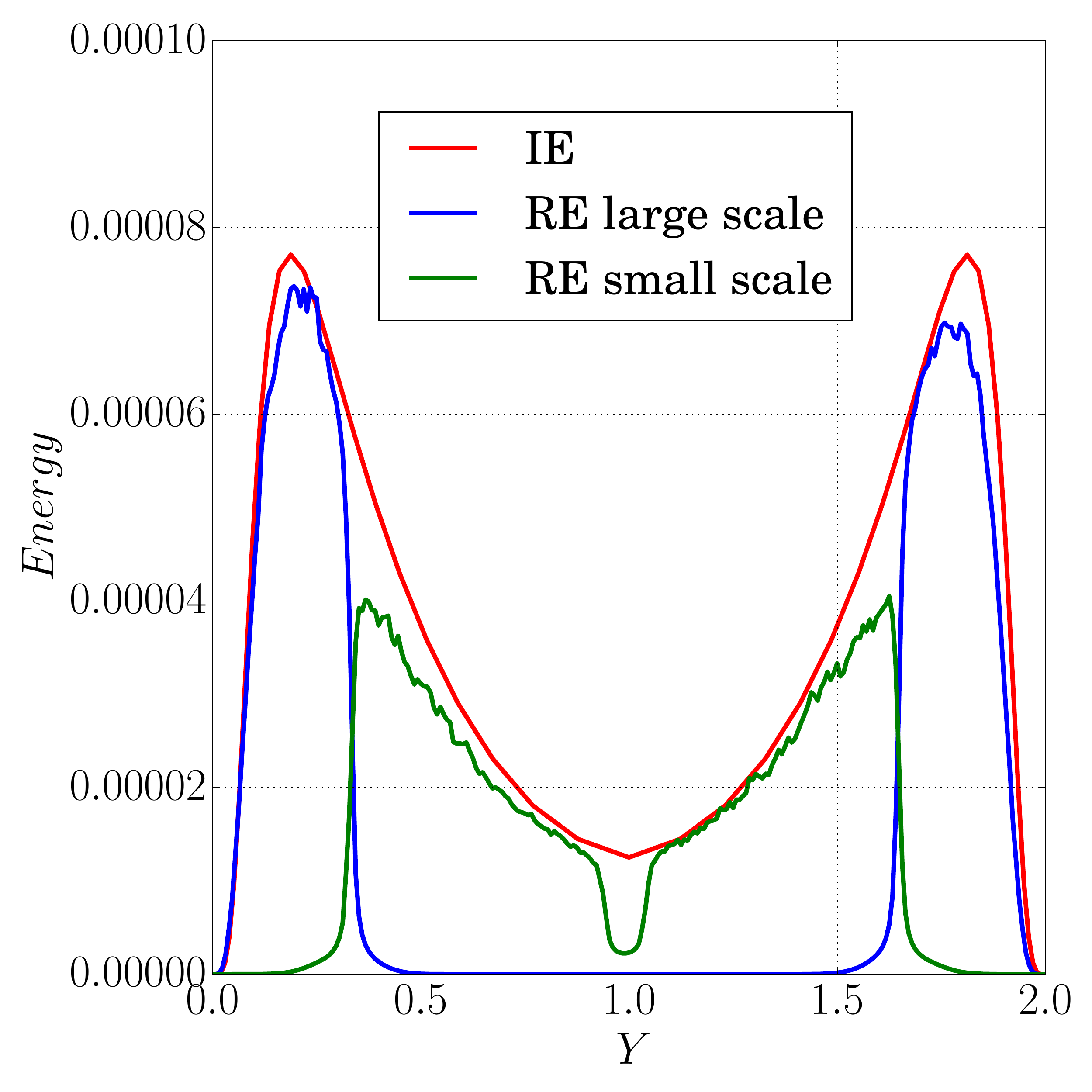}
\caption{$k=40.00$}
\label{fig:9sub7}
\end{subfigure}%
\begin{subfigure}{.5\textwidth}
\centering
\includegraphics[width=0.7\linewidth]{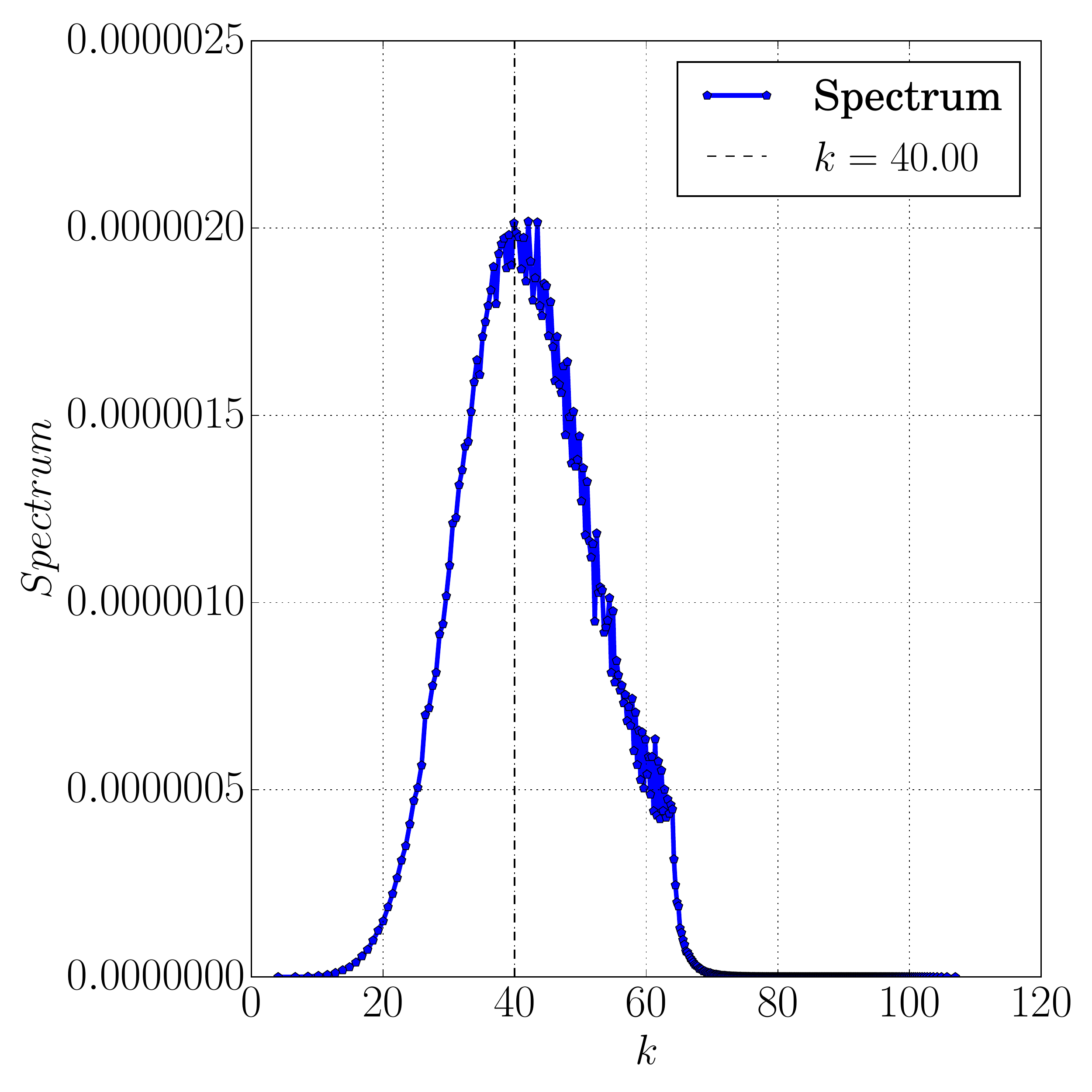}
\caption{$k=40.00$}
\label{fig:9sub8}
\end{subfigure}

\begin{subfigure}{.5\textwidth}
\centering
\includegraphics[width=0.7\linewidth]{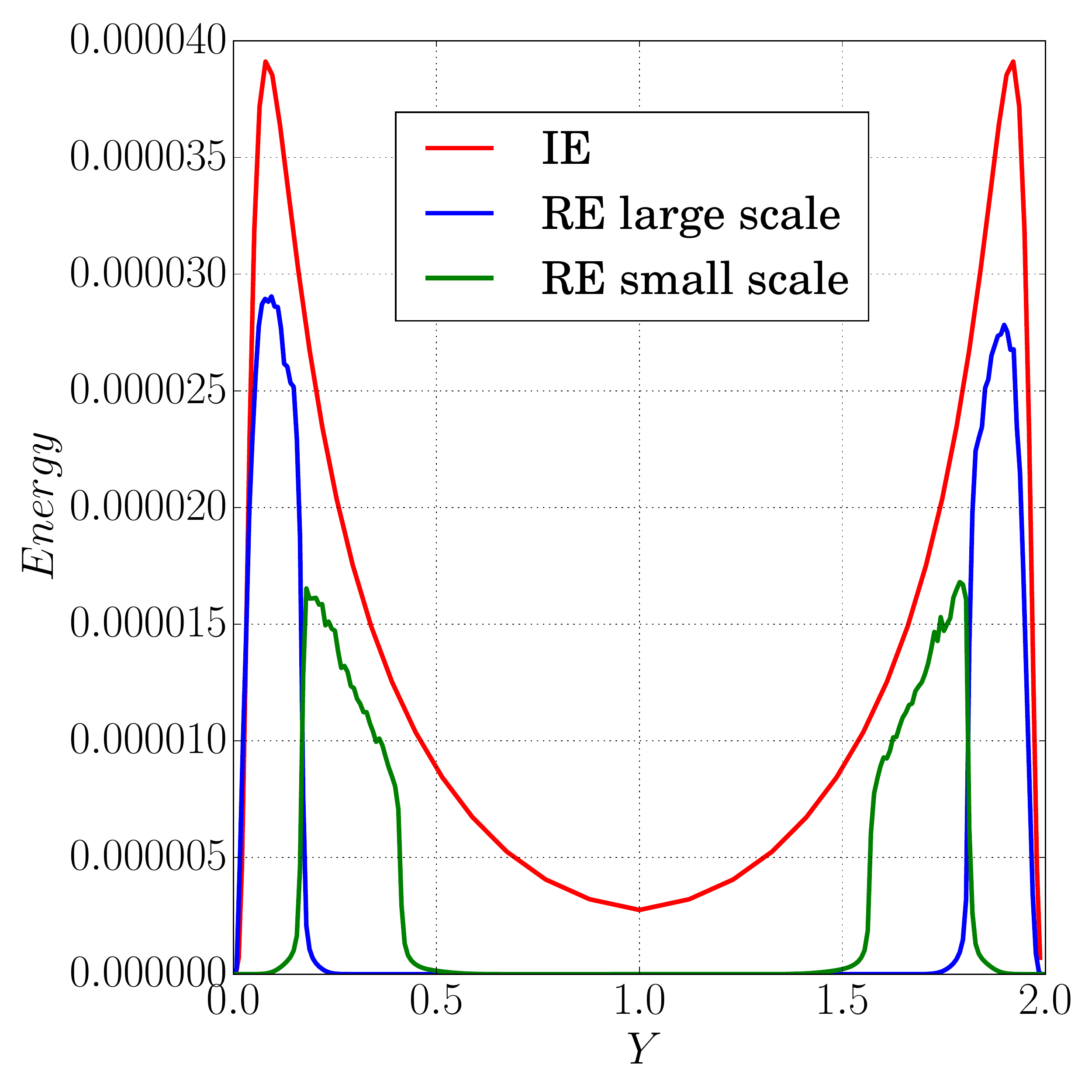}
\caption{$k=81.34$}
\label{fig:9sub9}
\end{subfigure}%
\begin{subfigure}{.5\textwidth}
\centering
\includegraphics[width=0.7\linewidth]{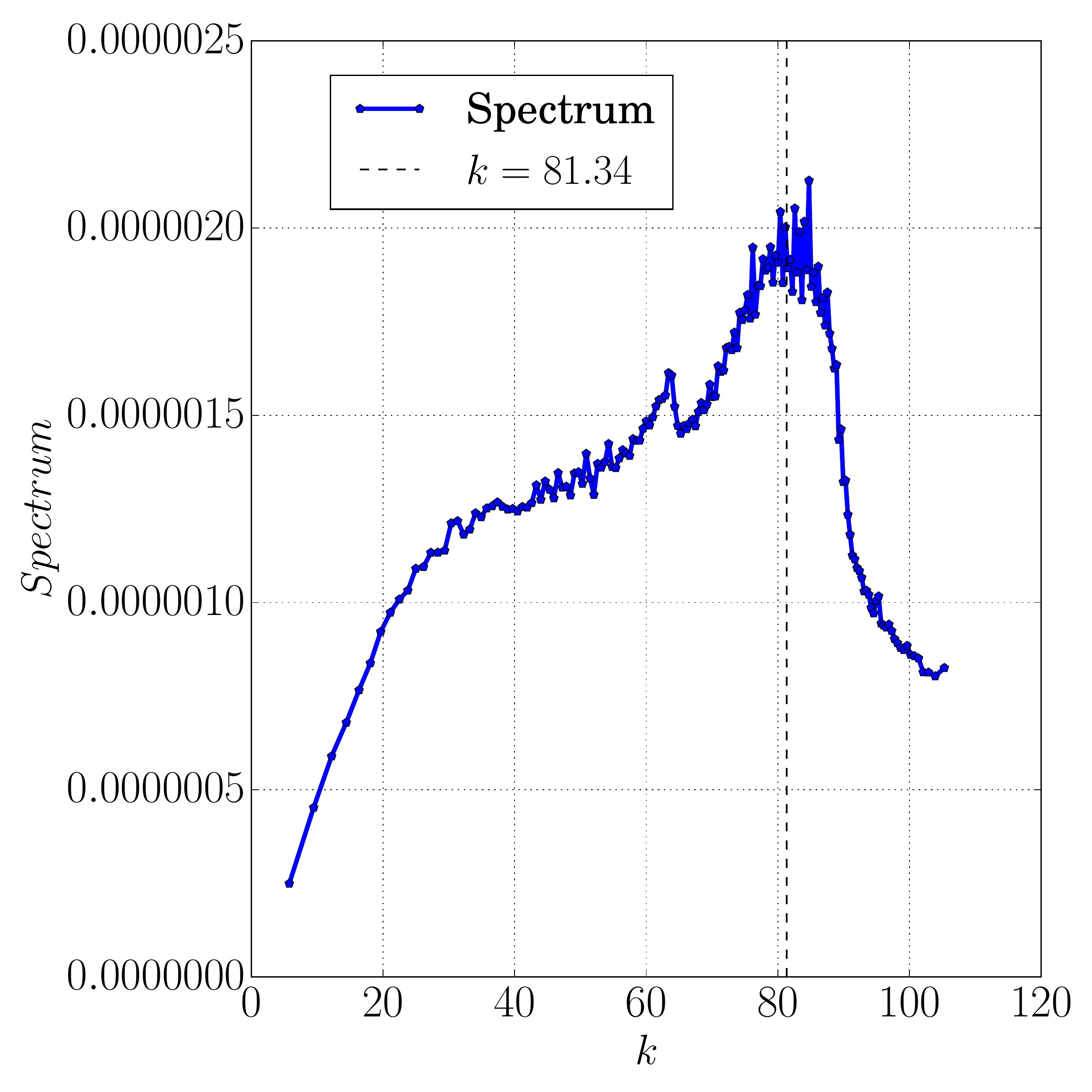}
\caption{$k=81.34$}
\label{fig:9sub10}
\end{subfigure}
\caption{Reconstruction of TKE at different wavenumber. Left column: spatial distribution of TKE at certain wavenumber. Red line: Input energy at certain $k$. Blue line: reconstructed energy of large scale at $k$. Green line:  reconstructed energy of small scale at $k$. Right column: spectral distribution of TKE. IE: Input Energy. RE: Reconstructed Energy. }
\label{fig:9}
\end{figure}

Fig.\ref{fig:9} shows TKE reconstruction at different wavenumbers. It can be observed that most of turbulent kinetic energy was fully reconstructed except for very high wavenumber case($k=81.34$). Also for each wavenumber, the energy of generated turbulence fluctuation concentrates at the neighbor of wavenumber of specific wavelet basis, which is the result of Eq.\ref{eq:7}, \ref{eq:8}. The reconstruction process is done in separated regions because of different time-advance schemes in Section 2.1.2. For $k=40.00$ and $k=81.34$, there is a region in which the reconstructed energy is $0$. It is because that TKE at this wavenumber only contributes to very small portion of total TKE at this region, thus is cutoff by the algorithm in Section 2.4. This portion of TKE can be accurately reconstructed, but the computation cost will increase significantly and the result does not change much. At very high wavenumber(Fig.\ref{fig:9sub9}, \ref{fig:9sub10}), although characteristic wavenumber of the wavelet mode is still below Nyquist wavenumber(in this case, $k_N$=120), some part of energy of the wavelet mode goes beyond $k_N$, which cannot be captured by the mesh in this case.

\begin{figure}[htbp]
\centering
\begin{subfigure}{.5\textwidth}
\centering
\includegraphics[width=0.7\linewidth]{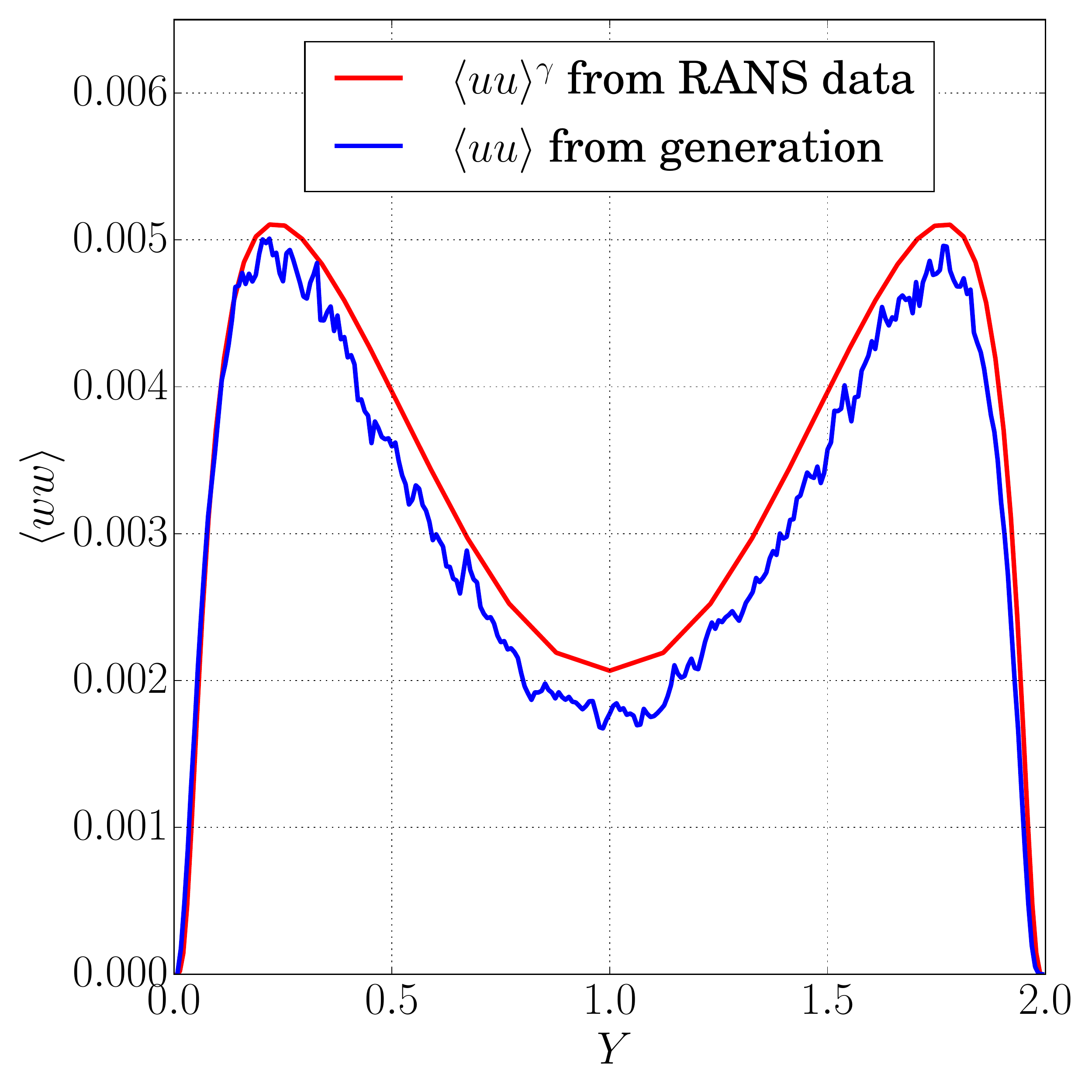}
\caption{$\langle uu \rangle$}
\label{fig:10sub1}
\end{subfigure}%
\begin{subfigure}{.5\textwidth}
\centering
\includegraphics[width=0.7\linewidth]{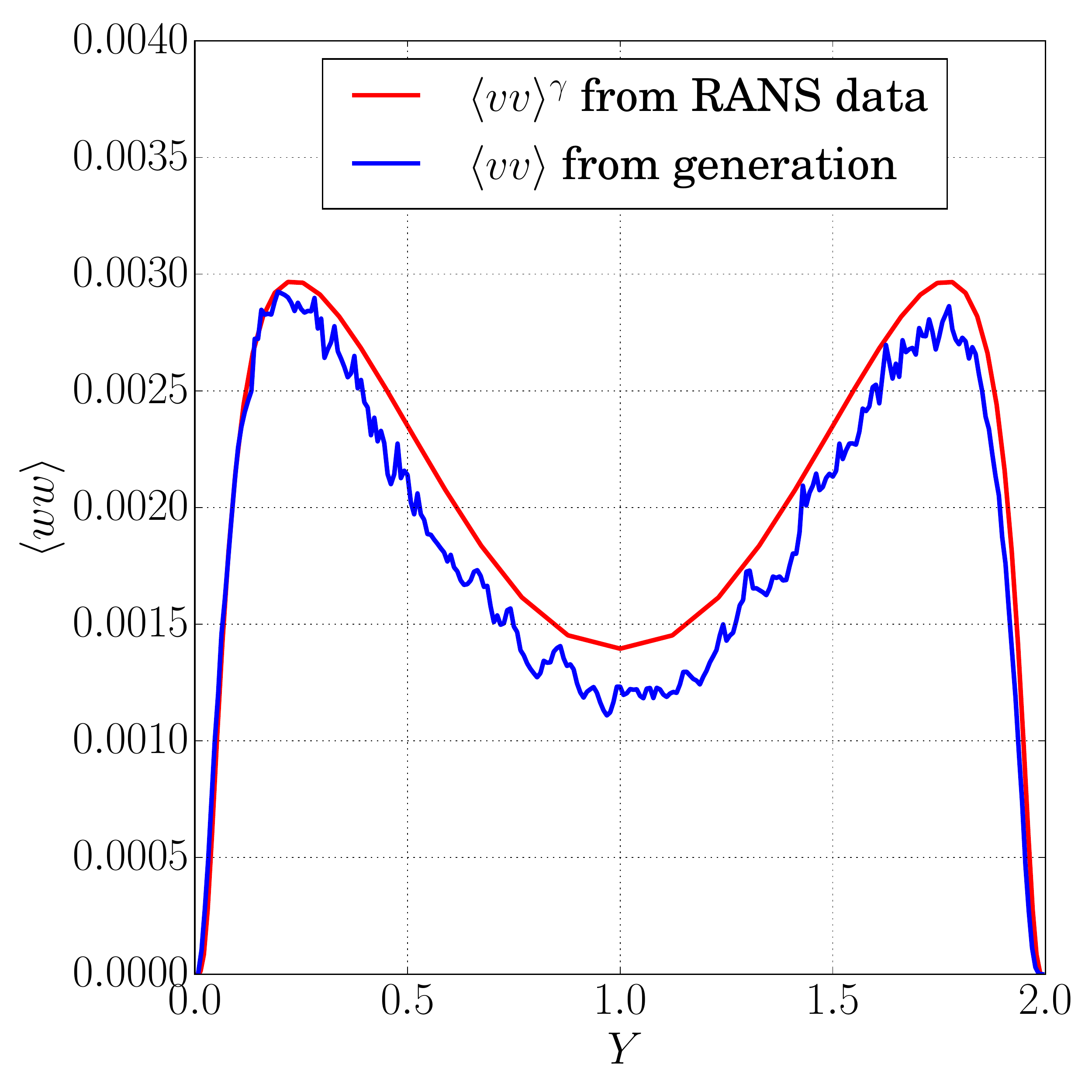}
\caption{$\langle vv \rangle$}
\label{fig:10sub2}
\end{subfigure}

\begin{subfigure}{.5\textwidth}
\centering
\includegraphics[width=0.7\linewidth]{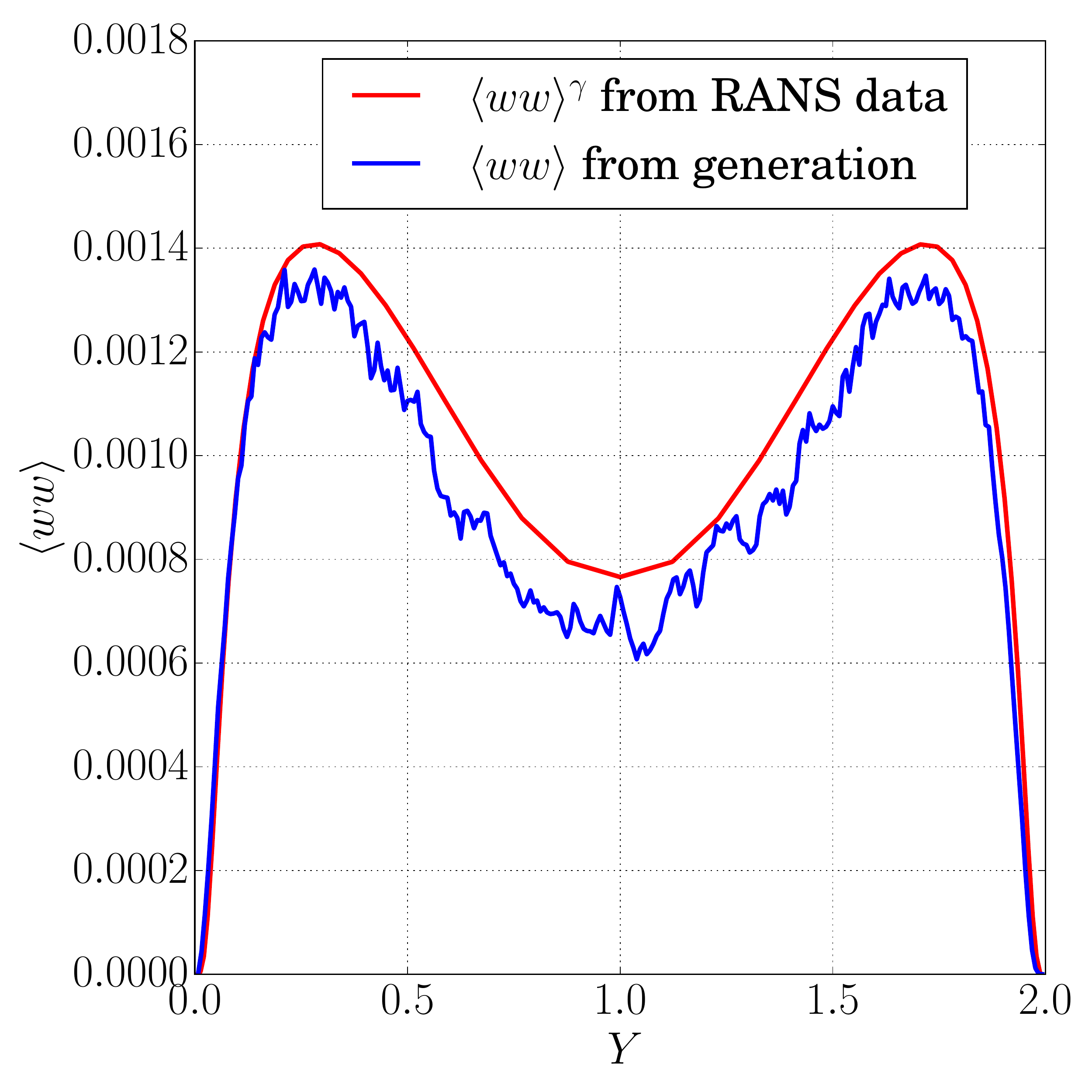}
\caption{$\langle ww \rangle$}
\label{fig:10sub3}
\end{subfigure}%
\begin{subfigure}{.5\textwidth}
\centering
\includegraphics[width=0.7\linewidth]{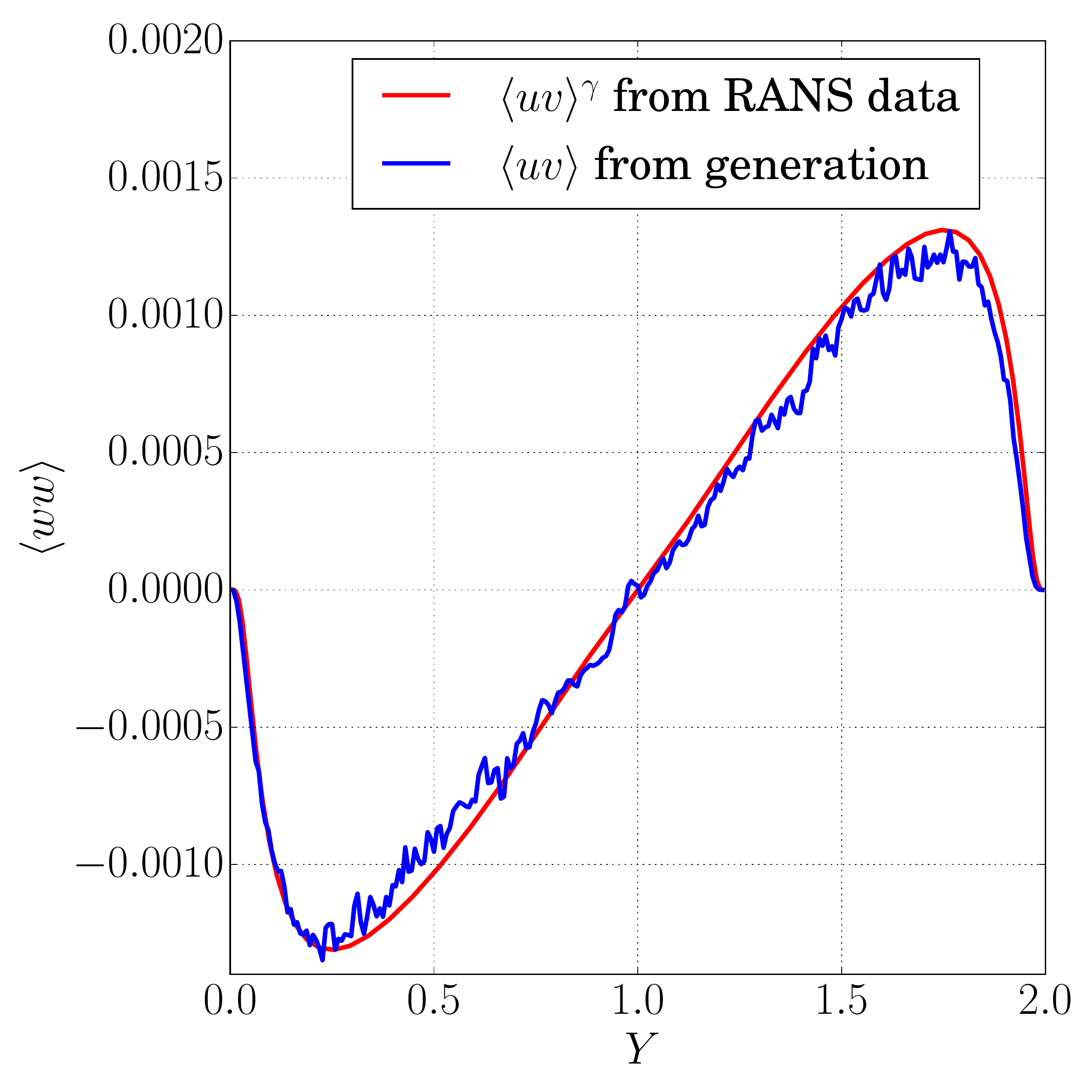}
\caption{$\langle uv \rangle$}
\label{fig:10sub4}
\end{subfigure}
\caption{Reconstruction of Reynolds stress distributions: The reconstructed Reynolds stress is compared with resolvable Reynolds stress from RANS data. Four main reconstructed Reynolds stresses all show good agreement with resolvable Reynolds stress}
\label{fig:10}
\end{figure}

Fig.\ref{fig:10} shows the comparison of Reynolds stress from RANS data and reconstructed Reynolds stress. It can be observed that reconstructed $\langle uu \rangle$, $\langle vv \rangle$, $\langle ww \rangle$, $\langle uv \rangle$ agree with RANS data. It should be noticed that RANS data used in Fig.\ref{fig:10} is slightly different from the data in Fig.\ref{fig:8}. Define reconstruction ratio as the follow:
$$
\gamma = \frac{\int_{0}^{k_{N}}E(k,y)\dif k}{k_t(y)}
$$
$\gamma$ represents the part of turbulent kinetic energy that could be resolved for given mesh. Then resolvable Reynolds stress $\langle \textrm{\textit{\textbf{u}}} \textrm{\textit{\textbf{u}}} \rangle ^{\gamma}$ is defined as follow:
$$
\langle \textrm{\textit{\textbf{u}}} \textrm{\textit{\textbf{u}}} \rangle ^{\gamma} = \gamma \langle \textrm{\textit{\textbf{u}}} \textrm{\textit{\textbf{u}}} \rangle
$$
$\langle \textrm{\textit{\textbf{u}}} \textrm{\textit{\textbf{u}}} \rangle ^{\gamma}$ represents the best approximation of Reynolds stress given a mesh of Nyquist wavenumber $k_{N}$. The reconstructed Reynolds stress in Fig.\ref{fig:10} shows good agreement with $\langle \textrm{\textit{\textbf{u}}} \textrm{\textit{\textbf{u}}} \rangle ^{\gamma}$ from RANS data. 

\begin{figure}[htbp]
\centering
\begin{subfigure}{.5\textwidth}
\centering
\includegraphics[width=1.0\linewidth]{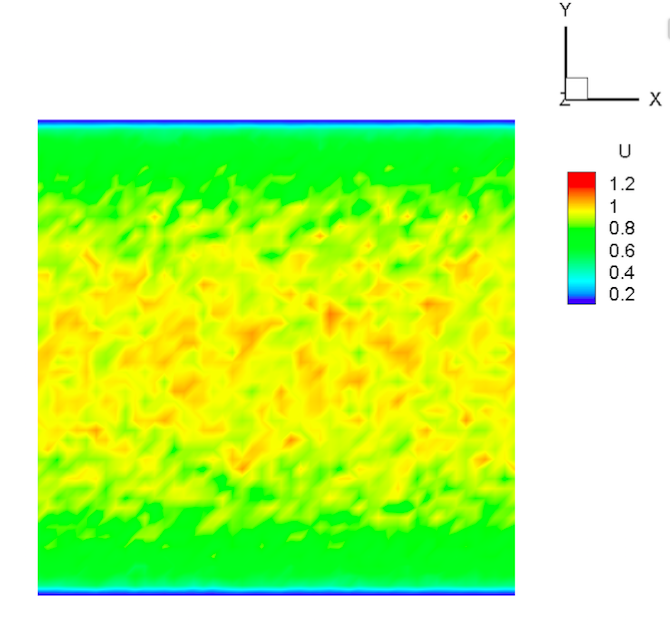}
\caption{Iso-contour of $U$}
\label{fig:11sub1}
\end{subfigure}%
\begin{subfigure}{.5\textwidth}
\centering
\includegraphics[width=1.0\linewidth]{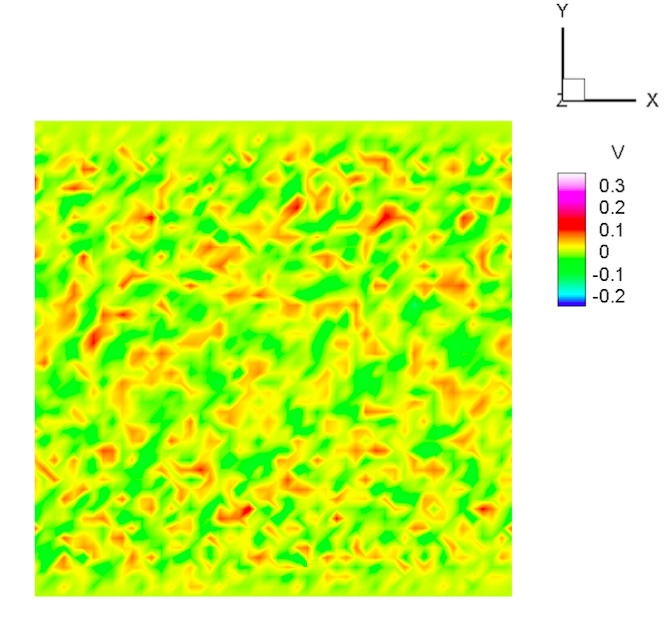}
\caption{Iso-contour of $V$}
\label{fig:11sub2}
\end{subfigure}
\begin{subfigure}{.5\textwidth}
\centering
\includegraphics[width=1.0\linewidth]{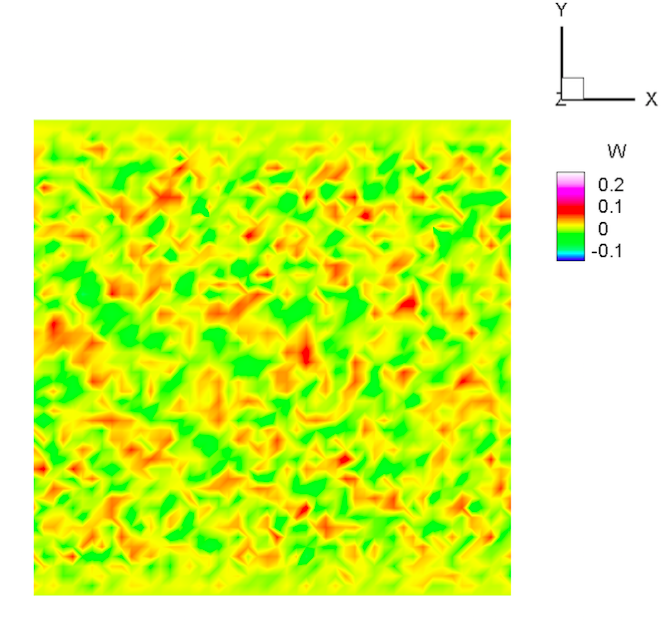}
\caption{Iso-contour of $W$}
\label{fig:11sub3}
\end{subfigure}
\caption{Iso-contour of different total velocity components of generated fully developed turbulent channel flow}
\label{fig:11}
\end{figure}       

Iso-contour of 3 total velocity components $U$, $V$, $W$ are shown in Fig.\ref{fig:11}. The reconstructed velocity field shows very realistic flow image. Large scale spatial structures are distributed near the centering line of channel. Near the wall turbulence structures get smaller and damping effect of the boundary starts to dominate. $V$ and $W$ components are close to 0, with some random fluctuation distributed in the cross section.

\section{Conclusion and Discussion}\label{Conclusion and Discussion}

In this paper, a new method of turbulence generation is proposed and evaluated in both homogeneous and inhomogeneous turbulence synthesis. Various properties of generated isotropic homogeneous turbulence shows good agreement with both input data and theoretical results, including spatial, spectral and frequency properties. Generated fully developed channel flow shows desired spectral and spatial characteristics for different wavenumber. Preservation of Reynolds Stress of this method is verified in both theoretical deduction and numerical simulation. 

Comparison of characteristics of Stochastic Wavelet Method and SRFM in homogeneous and inhomogeneous turbulence synthesis is listed in Tab.\ref{Tab:Comparison}. The number of modes of different wavenumbers used in turbulence synthesis with Stochastic Wavelet Model is far smaller than that with SRFM. Also, this new method could preserve local turbulence intensity as well as incompressibility, which could not be achieved with SRFM. Also, computational cost of Stochastic Wavelet Model could be largely reduced without much loss of turbulent kinetic energy, leading to far less computation cost than SRFM. 

This method exhibits great future in both scientific computing research and industrial application. Effective and low-cost inlet boundary generation is important for high-fidelity turbulence simulation(DNS, LES) and has become an important research topic recently. Also, interface between RANS region and LES region in zonal Detached Eddy Simulation also requires turbulence synthesis from RANS data. Stochastic Wavelet Method proposed in this paper provides a new approach to synthesize turbulence fluctuation field with desired spectral and statistical properties other than traditional SRFM. Also, this method could be further applied in Computer Graphics and movie industry to generate realistic fluid flow in animation with very low computation cost.

\begin{table}[]
\centering
\label{my-label}
\begin{tabular}{|c|c|c|c|c|c|c|}
\hline
\multirow{2}{*}{} & \multicolumn{2}{c|}{Homogeneous} & \multicolumn{4}{c|}{Inhomogeneous} \\ \cline{2-7} 
                  &    $\Re_L $     &  Mode Number        &  $\Re_{\tau}$   &   Reynolds Stress   &  Incompressibility   &  Mode Number   \\ \hline
SFRM                  &     723       &     5000      &  400   & Preserved    &   Not Preserved   &    500 - 5000     \\ \hline
SWM                  &      4206     &      10-20     &  $141900$   &  Preserved   & Preserved    &  15-30   \\ \hline
\end{tabular}
\caption{Comparison of SFRM and Stochastic Wavelet Model simulation results}
\label{Tab:Comparison}

\end{table}

\section{Appendix}\label{Appendix}
Consider turbulence field in a finite domain $\Omega$ of size $\vert \Omega \vert$. To prove the preservation of Reynolds stress tensor in scheme Eq.\ref{eq:1}, \ref{eq:2}, \ref{eq:3}, first to prove $\textrm{\textit{\textbf{v}}}= \nabla \times \textrm{\textit{\textbf{M}}}$ satisfies:
$$
\langle \textrm{\textit{\textbf{v}}} \textrm{\textit{\textbf{v}}} \rangle = \underline{\underline{I}}
$$
Proof:
$$
\textrm{\textit{\textbf{v}}} = \nabla \times \sum_{\vert \textrm{\textit{\textbf{k}}} \vert \in K}
       \sum^{N_i}_{\textrm{\textit{\textbf{x}}}_p}
       q_{\textrm{\textit{\textbf{x}}}_p,\textrm{\textit{\textbf{k}}}}\textrm{\textit{\textbf{O}}}_{\textrm{\textit{\textbf{x}}}_p,\textrm{\textit{\textbf{k}}}}(\mathbf{\omega}_{\textrm{\textit{\textbf{x}}}_p,\textrm{\textit{\textbf{k}}}} \Psi_{\textrm{\textit{\textbf{k}}}}(\textrm{\textit{\textbf{x}}}-\textrm{\textit{\textbf{x}}}_p))
$$

Curl is linear operator and is invariant under rotation:

$$
\textrm{\textit{\textbf{v}}} =  \sum_{\vert \textrm{\textit{\textbf{k}}} \vert \in K}
       \sum^{N_i}_{\textrm{\textit{\textbf{x}}}_p}
       q_{\textrm{\textit{\textbf{x}}}_p,\textrm{\textit{\textbf{k}}}}\textrm{\textit{\textbf{O}}}_{\textrm{\textit{\textbf{x}}}_p,\textrm{\textit{\textbf{k}}}}(\nabla \times(\mathbf{\omega}_{\textrm{\textit{\textbf{x}}}_p,\textrm{\textit{\textbf{k}}}} \Psi_{\textrm{\textit{\textbf{k}}}}(\textrm{\textit{\textbf{x}}}-\textrm{\textit{\textbf{x}}}_p)))
$$
Rewrite into index form:
$$
v_{i}=\sum_{\vert \textrm{\textit{\textbf{k}}} \vert \in K}
       \sum^{N_i}_{\textrm{\textit{\textbf{x}}}_p} q_{\textrm{\textit{\textbf{x}}}_p,\textrm{\textit{\textbf{k}}}} (\textrm{\textit{\textbf{O}}}_{\textrm{\textit{\textbf{x}}}_p,\textrm{\textit{\textbf{k}}}})_{il}\epsilon_{lmn}\partial_{m}(\mathbf{\omega}_{\textrm{\textit{\textbf{x}}}_p,\textrm{\textit{\textbf{k}}}} \Psi_{\textrm{\textit{\textbf{k}}}}(\textrm{\textit{\textbf{x}}}-\textrm{\textit{\textbf{x}}}_p))_{n}
$$
$$
= \sum_{\vert \textrm{\textit{\textbf{k}}} \vert \in K}
       \sum^{N_i}_{\textrm{\textit{\textbf{x}}}_p} q_{\textrm{\textit{\textbf{x}}}_p,\textrm{\textit{\textbf{k}}}} (\textrm{\textit{\textbf{O}}}_{\textrm{\textit{\textbf{x}}}_p,\textrm{\textit{\textbf{k}}}})_{il}\epsilon_{lmn}(\mathbf{\omega}_{\textrm{\textit{\textbf{x}}}_p,\textrm{\textit{\textbf{k}}}})_{n}\partial_{m}\Psi_{\textrm{\textit{\textbf{k}}}}(\textrm{\textit{\textbf{x}}}-\textrm{\textit{\textbf{x}}}_p)
$$
$$
\langle v_iv_j \rangle = 
\langle \sum_{\vert \textrm{\textit{\textbf{k}}}_1 \vert \in K}
       \sum^{N_i}_{\textrm{\textit{\textbf{x}}}_{p1}} 
       \sum_{\vert \textrm{\textit{\textbf{k}}}_2 \vert \in K}
       \sum^{N_i}_{\textrm{\textit{\textbf{x}}}_{p2}}
       q_{\textrm{\textit{\textbf{x}}}_{p1},\textrm{\textit{\textbf{k}}}_1} 
       q_{\textrm{\textit{\textbf{x}}}_{p2},\textrm{\textit{\textbf{k}}}_2} 
       (\textrm{\textit{\textbf{O}}}_{\textrm{\textit{\textbf{x}}}_{p1},\textrm{\textit{\textbf{k}}}_1})_{il}
       (\textrm{\textit{\textbf{O}}}_{\textrm{\textit{\textbf{x}}}_{p2},\textrm{\textit{\textbf{k}}}_2})_{jr}
      $$
      $$
       \epsilon_{lmn}
        \epsilon_{rst}
       (\mathbf{\omega}_{\textrm{\textit{\textbf{x}}}_{p1},\textrm{\textit{\textbf{k}}}_1})_{n}
       (\mathbf{\omega}_{\textrm{\textit{\textbf{x}}}_{p2},\textrm{\textit{\textbf{k}}}_2})_{t}
       \partial_{m}\Psi_{\textrm{\textit{\textbf{k}}}_1}(\textrm{\textit{\textbf{x}}}-\textrm{\textit{\textbf{x}}}_{p1})
       \partial_{s}\Psi_{\textrm{\textit{\textbf{k}}}_2}(\textrm{\textit{\textbf{x}}}-\textrm{\textit{\textbf{x}}}_{p2})\rangle
      $$
      $$
      =
 \sum_{\vert \textrm{\textit{\textbf{k}}}_1 \vert \in K}
       \sum^{N_i}_{\textrm{\textit{\textbf{x}}}_{p1}} 
       \sum_{\vert \textrm{\textit{\textbf{k}}}_2 \vert \in K}
       \sum^{N_i}_{\textrm{\textit{\textbf{x}}}_{p2}}
       q_{\textrm{\textit{\textbf{x}}}_{p1},\textrm{\textit{\textbf{k}}}_1} 
       q_{\textrm{\textit{\textbf{x}}}_{p2},\textrm{\textit{\textbf{k}}}_2} 
       \langle(\textrm{\textit{\textbf{O}}}_{\textrm{\textit{\textbf{x}}}_{p1},\textrm{\textit{\textbf{k}}}_1})_{il}
       (\textrm{\textit{\textbf{O}}}_{\textrm{\textit{\textbf{x}}}_{p2},\textrm{\textit{\textbf{k}}}_2})_{jr}\rangle
      $$
      $$
      \epsilon_{lmn}
        \epsilon_{rst}
       \langle(\mathbf{\omega}_{\textrm{\textit{\textbf{x}}}_{p1},\textrm{\textit{\textbf{k}}}_1})_{n}
       (\mathbf{\omega}_{\textrm{\textit{\textbf{x}}}_{p2},\textrm{\textit{\textbf{k}}}_2})_{t}\rangle
       \langle\partial_{m}\Psi_{\textrm{\textit{\textbf{k}}}_1}(\textrm{\textit{\textbf{x}}}-\textrm{\textit{\textbf{x}}}_{p1})
       \partial_{s}\Psi_{\textrm{\textit{\textbf{k}}}_2}(\textrm{\textit{\textbf{x}}}-\textrm{\textit{\textbf{x}}}_{p2})\rangle
      $$
      For $p_1 \neq p_2$ or $ k_1 \neq k_2$, $\langle(\mathbf{\omega}_{\textrm{\textit{\textbf{x}}}_{p1},\textrm{\textit{\textbf{k}}}_1})_{n}
       (\mathbf{\omega}_{\textrm{\textit{\textbf{x}}}_{p2},\textrm{\textit{\textbf{k}}}_2})_{t}\rangle=0$. Thereforth:
       $$
       \langle v_iv_j \rangle = \sum_{\vert \textrm{\textit{\textbf{k}}}_ \vert \in K}
       \sum^{N_i}_{\textrm{\textit{\textbf{x}}}_{p}} 
       q_{\textrm{\textit{\textbf{x}}}_{p},\textrm{\textit{\textbf{k}}}} ^2
       \langle(\textrm{\textit{\textbf{O}}}_{\textrm{\textit{\textbf{x}}}_{p},\textrm{\textit{\textbf{k}}}})_{il}
       (\textrm{\textit{\textbf{O}}}_{\textrm{\textit{\textbf{x}}}_{p},\textrm{\textit{\textbf{k}}}})_{jr}\rangle
       \epsilon_{lmn}
        \epsilon_{rst}
      $$
      $$
       \langle(\mathbf{\omega}_{\textrm{\textit{\textbf{x}}}_{p},\textrm{\textit{\textbf{k}}}})_{n}
       (\mathbf{\omega}_{\textrm{\textit{\textbf{x}}}_{p},\textrm{\textit{\textbf{k}}}})_{t}\rangle
       \langle\partial_{m}\Psi_{\textrm{\textit{\textbf{k}}}}(\textrm{\textit{\textbf{x}}}-\textrm{\textit{\textbf{x}}}_{p})
       \partial_{s}\Psi_{\textrm{\textit{\textbf{k}}}}(\textrm{\textit{\textbf{x}}}-\textrm{\textit{\textbf{x}}}_{p})\rangle
      $$
      $\textrm{\textit{\textbf{x}}}_{p}$ is uniformly distributed in the flow domain, thus the following hold:
      $$
      \langle\partial_{m}\Psi_{\textrm{\textit{\textbf{k}}}}(\textrm{\textit{\textbf{x}}}-\textrm{\textit{\textbf{x}}}_{p})
       \partial_{s}\Psi_{\textrm{\textit{\textbf{k}}}}(\textrm{\textit{\textbf{x}}}-\textrm{\textit{\textbf{x}}}_{p})\rangle = \langle\int_{\Omega}\partial_m\Psi_{\textrm{\textit{\textbf{k}}}}\partial_s\Psi_{\textrm{\textit{\textbf{k}}}}\dif \textrm{\textit{\textbf{x}}}\rangle
      $$
      From construction process in Eq.\ref{eq:7}, $\partial_{s}\Psi_{\textrm{\textit{\textbf{k}}}}$ is symmetric along three axis, thus:
      $$
      \langle\int_{\Omega}\partial_m\Psi_{\textrm{\textit{\textbf{k}}}}\partial_s\Psi_{\textrm{\textit{\textbf{k}}}}\dif \textrm{\textit{\textbf{x}}}\rangle = \frac{c_k }{N_i}\delta_{ms}
      $$
      where $c_k$ is defined in Eq.\ref{eq:11}. Also notice $\langle(\mathbf{\omega}_{\textrm{\textit{\textbf{x}}}_{p},\textrm{\textit{\textbf{k}}}})_{n}
       (\mathbf{\omega}_{\textrm{\textit{\textbf{x}}}_{p},\textrm{\textit{\textbf{k}}}})_{t}\rangle=\delta_{nt}$ Plug into the expression of $\langle v_iv_j \rangle$:
      $$
      \langle v_iv_j \rangle = \sum_{\vert \textrm{\textit{\textbf{k}}} \vert \in K}
       \sum^{N_i}_{\textrm{\textit{\textbf{x}}}_{p}} 
       q_{\textrm{\textit{\textbf{x}}}_{p},\textrm{\textit{\textbf{k}}}} ^2
       \langle(\textrm{\textit{\textbf{O}}}_{\textrm{\textit{\textbf{x}}}_{p},\textrm{\textit{\textbf{k}}}})_{il}
       (\textrm{\textit{\textbf{O}}}_{\textrm{\textit{\textbf{x}}}_{p},\textrm{\textit{\textbf{k}}}})_{jr}\rangle
       \epsilon_{lmn}
        \epsilon_{rst}
        \delta_{nt}\delta_{ms}\frac{c_k }{N_i}
      $$
      $$
      = \sum_{\vert \textrm{\textit{\textbf{k}}} \vert \in K}
       \sum^{N_i}_{\textrm{\textit{\textbf{x}}}_{p}} 
       2 q_{\textrm{\textit{\textbf{x}}}_{p},\textrm{\textit{\textbf{k}}}} ^2
       \langle(\textrm{\textit{\textbf{O}}}_{\textrm{\textit{\textbf{x}}}_{p},\textrm{\textit{\textbf{k}}}})_{il}
       (\textrm{\textit{\textbf{O}}}_{\textrm{\textit{\textbf{x}}}_{p},\textrm{\textit{\textbf{k}}}})_{jr}\rangle
       \delta_{lr} \frac{ c_k}{N_i}
      $$
      $$
      = \sum_{\vert \textrm{\textit{\textbf{k}}} \vert \in K}
       \sum^{N_i}_{\textrm{\textit{\textbf{x}}}_{p}} 
       2 q_{\textrm{\textit{\textbf{x}}}_{p},\textrm{\textit{\textbf{k}}}} ^2
       \langle(\textrm{\textit{\textbf{O}}}_{\textrm{\textit{\textbf{x}}}_{p},\textrm{\textit{\textbf{k}}}})_{il}
       (\textrm{\textit{\textbf{O}}}_{\textrm{\textit{\textbf{x}}}_{p},\textrm{\textit{\textbf{k}}}})_{jl}\rangle
       \frac{c_k}{N_i}
      $$
      $$
      = \sum_{\vert \textrm{\textit{\textbf{k}}} \vert \in K}
       \sum^{N_i}_{\textrm{\textit{\textbf{x}}}_{p}} 
       2 q_{\textrm{\textit{\textbf{x}}}_{p},\textrm{\textit{\textbf{k}}}} ^2
       \delta_{ij}
       \frac{ c_k}{N_i}
      $$
      $$
      =\delta_{ij} \sum_{\vert \textrm{\textit{\textbf{k}}} \vert \in K}
        \frac{E(l)\Delta l}{2k_t } \to \delta_{ij} 
        $$
        as $\Delta l \to 0, l_{max} \to \infty, l_{min} \to 0  $. Thus for enough large $K$, the following hold:
        $$
        \langle \textrm{\textit{\textbf{v}}} \textrm{\textit{\textbf{v}}} \rangle = \textrm{\textit{\textbf{I}}}
        $$
      Then:
      $$
      \langle \textrm{\textit{\textbf{u}}} \textrm{\textit{\textbf{u}}} \rangle = \langle  (\textrm{\textit{\textbf{A}}} \textrm{\textit{\textbf{v}}})(\textrm{\textit{\textbf{A}}} \textrm{\textit{\textbf{v}}})\rangle
      $$
      $$
       = \textrm{\textit{\textbf{A}}}(\textrm{\textit{\textbf{v}}} \textrm{\textit{\textbf{v}}}) \textrm{\textit{\textbf{A}}}^{T}
      $$
      $$
       = \textrm{\textit{\textbf{A}}}\textrm{\textit{\textbf{A}}}^{T} = \textrm{\textit{\textbf{R}}}
      $$

\newpage

\bibliographystyle{unsrt}
\bibliography{citation}
\end{document}